\documentclass[10pt]{amsart}

\usepackage{amssymb}
\usepackage{amsmath}
\usepackage{amsthm}
\usepackage{amssymb}

\begin{document}

\swapnumbers
\theoremstyle{plain}
\newtheorem{theorem}[subsection]{Theorem}
\newtheorem{assumption}[subsection]{Assumption}
\newtheorem*{theorem*}{Theorem}
\newtheorem{proposition}[subsection]{Proposition}
\newtheorem{lemma}[subsection]{Lemma}
\newtheorem*{lemma*}{Lemma}
\newtheorem{corollary}[subsection]{Corollary}
\newtheorem*{conjecture}{Conjecture}
\theoremstyle{definition}
\newtheorem*{definition}{Definition}
\theoremstyle{remark}
\newtheorem*{remark}{Remark}
\newtheorem*{remarks}{Remarks}
\newtheorem{example}[subsubsection]{Example}
\newtheorem{claim}[subsubsection]{Claim}
\newtheorem*{acknowledgements}{Acknowledgements}

\renewcommand{\descriptionlabel}[1]%
         {\hspace{\labelsep}\normalfont{#1}}

\newcommand{\0}{{(0)}}
\renewcommand{\AA}{\mathbf{A}}
\newcommand{\BB}{\mathbf{B}}
\newcommand{\CC}{\mathbb{C}}
\newcommand{\DD}{\mathbf{D}}
\newcommand{\FF}{\mathbb{F}}
\newcommand{\Ff}{\mathcal{F}}
\newcommand{\GG}{\mathbf{G}}
\newcommand{\HH}{\mathbf{H}}
\newcommand{\LL}{\mathbf{L}}
\newcommand{\MM}{\mathbf{M}}
\newcommand{\NN}{\mathbf{N}}
\newcommand{\PP}{\mathbf{P}}
\renewcommand{\SS}{\mathbf{S}}
\newcommand{\TT}{\mathbf{T}}
\newcommand{\UU}{\mathbf{U}}
\newcommand{\VV}{\mathbf{V}}
\newcommand{\XX}{\mathbf{X}}
\newcommand{\YY}{\mathbf{Y}}
\newcommand{\ZZ}{\mathbf{Z}}
\newcommand{\QQ}{\mathbb{Q}}

\newcommand{\ur}{\textrm{ur}}
\newcommand{\rr}{\mathcal{r}}

\newcommand{\cind}{\operatorname{c-Ind}}
\newcommand{\Ind}{\operatorname{Ind}}
\newcommand{\Hom}{\operatorname{Hom}}
\newcommand{\End}{\operatorname{End}}
\newcommand{\Aut}{\operatorname{Aut}}
\newcommand{\varchi}{\mathcal{X}}
\newcommand{\Gm}{\mathbb{G}_m}
\newcommand{\Ga}{\mathbb{G}_a}
\newcommand{\supp}{\operatorname{supp}}
\newcommand{\GGL}{\operatorname{\mathbf{GL}}}
\newcommand{\GL}{\operatorname{GL}}
\newcommand{\PPGL}{\operatorname{\mathbf{\mathbf{PGL}}}}
\newcommand{\PGL}{\operatorname{PGL}}
\newcommand{\SSL}{\operatorname{\mathbf{SL}}}
\newcommand{\SL}{\operatorname{SL}}
\newcommand{\Sp}{\operatorname{Sp}}
\newcommand{\Ssp}{\operatorname{\mathbf{Sp}}}
\newcommand{\SO}{{\operatorname{SO}}}
\newcommand{\geom}{\operatorname{geom}}
\newcommand{\Gal}{\operatorname{Gal}}
\newcommand{\Fr}{{\operatorname{(Fr)}}}
\newcommand{\tr}{\operatorname{tr}}
\newcommand{\Adm}{\mathit{Adm}}
\newcommand{\spec}{\operatorname{spec}}
\newcommand{\Vol}{\operatorname{Vol}}
\newcommand{\diag}{{\operatorname{diag}}}
\newcommand{\ev}{\operatorname{ev}}
\newcommand{\meas}{\operatorname{\meas}}

\newcommand{\nc}{\newcommand}
\nc{\cali}{\mathcal} 
\nc{\la}{\langle} \nc{\ra}{\rangle}
 \nc{\CA}{\cali A}
 \nc{\CBB}{\cali B}
\nc{\CDD}{\cali D}
\nc{\CE}{\cali E}
\nc{\CF}{\cal F} \nc{\CG}{\cal
G} \nc{\CH}{\cali H} \nc{\CI}{\cal I} \nc{\CJ}{\cal J}
\nc{\CK}{\cal K} \nc{\CL}{\cal L} \nc{\CM}{\cal M} \nc{\CN}{\cal
N} \nc{\CO}{\cal O} \nc{\CP}{\cal P} \nc{\CQ}{\cal Q}
\nc{\CR}{\cal R} \nc{\CS}{\cal S} \nc{\CT}{\cal T} \nc{\CU}{\cal
U} \nc{\CV}{\cal V} \nc{\CW}{\cal W} \nc{\CZ}{\cal Z}


\nc{\fa}{\mathfrak a} \nc{\fg}{\mathfrak g} \nc{\fk}{\mathfrak k}
\nc{\fh}{\mathfrak h} \nc{\fm}{\mathfrak m} \nc{\fn}{\mathfrak n}
\nc{\fA}{\mathfrak A} \nc{\fC}{\mathfrak C} \nc{\fI}{\mathfrak I}
\nc{\fL}{\mathfrak L} \nc{\fS}{\mathfrak S}


\nc{\nen}{\newenvironment} \nc{\ol}{\overline}
\nc{\ul}{\underline} \nc{\lra}{\longrightarrow}
\nc{\lla}{\longleftarrow} \nc{\Lra}{\Longrightarrow}
\nc{\Lla}{\Longleftarrow} \nc{\Llra}{\Longleftrightarrow}
\nc{\hra}{\hookrightarrow} \nc{\iso}{\overset{\sim}{\lra}}

\makeatletter
\@addtoreset{equation}{section}
\makeatother
\numberwithin{equation}{section}
 \nc{\ba}{\mathbb A}
 \nc{\bq}{\mathbb Q}
 \nc{\br}{\mathbb R}
 \nc{\bz}{\mathbb Z}
 \nc{\bc}{\mathbb C}
 \nc{\bn}{\mathbb N}
 \nc{\ck}{\mathcal{K}}
 \nc{\G}{\Gamma}
 \nc{\sm}{\setminus}
 \nc{\sub}{\subset}
 \nc{\lm}{\lambda}
  \nc{\Lm}{\Lambda}
 \nc{\al}{\alpha}
 \nc{\bt}{\beta}
 \nc{\om}{\omega}
 \nc{\dl}{\delta}
 \nc{\g}{\gamma}
 \nc{\Dl}{\Delta}
 \nc{\Om}{\Omega}
 \nc{\s}{\sigma}
 \nc{\ro}{\rho}
 \nc{\te}{\theta}
 \nc{\SLR}{\operatorname{SL}_2(\br)}
 \nc{\GLR}{\operatorname{GL}_2(\br)}
 \nc{\PGLR}{\operatorname{PGL}_2(\br)}
 \nc{\PSLR}{\operatorname{PSL}_2(\br)}
 \nc{\SLC}{\operatorname{SL}_2(,\bc)}
 \nc{\POt}{\operatorname{PO}(2)}
 \nc{\PSOt}{\operatorname{PSO}(2)}
 \nc{\Ot}{\operatorname{O}(2)}
 \nc{\uH}{\mathbb H}
 \nc{\fD}{\mathcal{D}}
 \nc{\fE}{\mathcal{E}}
 \nc{\fO}{\mathcal{O}}
 \nc{\haf}{\frac{1}{2}}
 \nc{\qtr}{\frac{1}{4}}
 \nc{\shaf}{{\scriptstyle\frac{1}{2}}}
 \nc{\hlm}{{\scriptstyle\frac{\lambda}{2}}}

 \nc{\8}{\infty}
 \nc{\7}{{-\infty}}
 \nc{\inv}{^{-1}}
 \nc{\eps}{\varepsilon}
 \nc{\aG}{\mathbf{G}}
 \nc{\spn}{\operatorname{Span}}


\setcounter{tocdepth}{1}

\title[Casimir operator]{Geodesic restrictions for the Casimir operator}

\author{Andre Reznikov} \email{reznikov@math.biu.ac.il}

\address{Department of Mathematics\\
Bar-Ilan University \\
Ramat Gan, 52900 \\
Israel}

\begin{abstract} We consider the hyperbolic Casimir operator $C$ defined on the tangent sphere bundle $SY$ of a compact hyperbolic Riemann surface $Y$. We prove a nontrivial bound on the $L^2$-norm of the restriction of eigenfunctions of $C$ to certain natural hypersurfaces in $SY$. The result that we obtain goes beyond known (sharp) local bounds of L. H\"{o}rmander.
\end{abstract}


\thanks{Partially supported by the Veblen Fund at IAS, by a BSF grant, by the ISF  Center of Excellency grant 1691/10, and by the Minerva  Center at ENI.}

\maketitle

\section{Introduction}\label{secintro}
\subsection{Restriction problem} The purpose of this paper is to draw attention to the following well-known question in PDE.
Let $X$ be a compact $n$-dimensional smooth manifold, and let $C$ be a second order hyperbolic operator (i.e., an operator with the principal symbol $p$ of the signature $(n-1, 1)$). Let $O\subset X$ be a closed smooth hypersurface, which we assume to be non-characteristic for $C$.  Let $f\in C^\8(X)$ be a smooth function. What can be said about the norm of the restriction of $f$ to $O$ in terms of norms of $f$ and $C(f)$ on $X$?

The ``local" theory related to such questions goes back to works of L.~H\"{o}rmander (for example, see \cite[Theorem 25.3.11]{Ho}). One knows that the answer depends on the curvature of $O$ with respect to the pseudometric defined by the principal symbol $p$. Resulting bounds and various  extensions  are discussed at length in the paper by D. Tataru, \cite{Ta}.
In particular, for the smooth time-like hypersurface we have
\begin{equation}\label{00}
 ||f|_O||_{W^{3/4}_{loc}(O)}\leq A\cdot\left(||f||_{W^1_{loc}(X)}+||C(f)||_{L^2_{loc}(X)}\right)\ ,
  \end{equation}for some constant $A>0$ depending on the geometry of $X$, $O$ and $p$. Here we denote by $W^s$  the $L^2$-Sobolev norm of the order $s$ on $X$ (for example, associated to a choice of Riemannian  metric on $X$). This should be compared with the Sobolev embedding theorem giving the loss of $1/2$ of the derivative (i.e.,
  $||f|_O||_{W^{1/2}_{loc}(O)}\leq A'\cdot||f||_{W^1_{loc}(X)}$).  If in addition $O$ is curved, A. Greenleaf and A. Seeger \cite{GrS} obtained the improvement  $||f|_O||_{W^{5/6}_{loc}(O)}\leq A''\cdot\left(||f||_{W^1_{loc}(X)}+||C(f)||_{L^2_{loc}(X)}\right)$ of the bound \eqref{00}.  These bounds are optimal, as could be seen from the constant coefficients examples.

Since these results are {\it local}, one can not hope to improve these bounds for general functions using global geometry of $X$, $O$ and $p$. However, we want to present an example where the improvement is possible for functions having dominating first term on the right in the bound \eqref{00}. Namely, we will construct an example where
 \begin{equation}\label{01}
 ||f|_O||_{L^2(O)}\leq A_{a,\eps}\cdot||f||_{W^\eps(X)}\ ,
  \end{equation}for any $\eps>0$, and for any $f\in C^\8(X)$ satisfying $||C(f)||_{W^\eps(X)}\leq a$ for a given $a>0$. In particular, for eigenfunctions of $C$ with bounded eigenvalue,  there is essentially no loss of smoothness in taking the restriction (at least to some special hypersurfaces). This raises the question if such a phenomenon persists more generally. Namely, one would like to see if it is possible to impose certain conditions on the geometry of $X$, $O$ and $p$ which would imply that
\begin{equation}\label{0}
||f|_O||_{L^2(O)}\leq A_{X,O,p,\eps}\cdot\left(||f||_{W^\eps(X)}+||C(f)||_{W^\eps(X)}\right)\
\end{equation} for {\it eigenfunctions} of $C$.

As we mentioned above, one can not hope for a local explanation for a bound of the type \eqref{0}. In fact, one can easily construct an example  (based on spherical harmonics, as usual) where nothing close to the bound \eqref{0}  holds for a compact $X$ and $O$ (see Example \ref{S2} below). Nevertheless, an intuition from recent works on Quantum Chaos suggests that there are some natural situations where one might expect some version of \eqref{0} to hold. This paper proposes one such an example which is based on the theory of hyperbolic Riemann surfaces (or more appropriately, on representation theory of $\PGLR$). We show that in that particular example the bound \eqref{0} holds for the ``low energy" spectrum, and for special hypersurfaces $O$ (e.g., associated to closed geodesic circles on the Riemann surface).

Our proof (and the example) comes from  representation theory, and is a mixture of results and techniques from \cite{1} and \cite{2}.

\subsection{Casimir operator} Here we describe our construction. Namely, we construct a $3$-dimensional compact manifold $X$ ($S^1$-fibered over a compact hyperbolic Riemann surface $Y$), and a nondegenerate  second order hyperbolic operator $C$  on $X$ (the Casimir operator) such that for any immersed two-dimensional tori $O\subset X$ which is fibered over a {\it closed} geodesic or a {\it closed} geodesic circle on $Y$, the bound \eqref{01} holds.

First, let us recall the well-known setup  (due to I. Gelfand and S. Fomin, \cite{G6}) which allows one to apply the representation theory to some special PDE operators.
\subsubsection{Hyperbolic geometry: Riemann surfaces and Geodesics}
Consider  a compact Riemann surface $(Y,g_2)$ endowed with the constant negative curvature metric $g_2$,
the corresponding volume element $dv_2$, and let $\Delta$ be the  corresponding ({\it nonnegative}) Laplace operator on $Y$. Let $SY$ be the tangent sphere bundle over $Y$. It is more convenient for our purposes to work with the oriented orthonormal frame bundle  $RY$ over $Y$. The $3$-dimensional manifold $RY$ has two connected components for the oriented Riemann surface $Y$, and there is the projection map $RY\to SY$ taking the first vector in the frame.

The central object in our method is the group $G=\PGLR$ (one can work with the group $\operatorname{SL}_2^\pm(\br)$ of matrices $g$ with $\det(g)=\pm 1$ instead). This group naturally appears in geometry as the group of isometries of the hyperbolic plane $\CH$ not necessarily preserving orientation. We fix the standard maximal compact subgroup $K=\POt\subset G$, and denote by $K_0=\PSOt\subset \operatorname{PGL}_2^+(\br)$ its maximal connected  subgroup ($K_0\simeq S^1$). We have the isomorphism $\CH\simeq G/K$. The uniformization theorem implies that there exists a lattice  $\Gamma\subset G$ (i.e., a discrete co-compact subgroup $\G\simeq\pi_1(Y)$) such that $Y=\G\setminus\CH$.  We set  $X=\G\setminus G$.  This will be our manifold. It is important to note that $X$ is a homogeneous manifold for $G$ acting on the right. We have the isomorphism $X\simeq RY$. We will assume for simplicity that $\G\subset \PSLR$ (i.e., we distinguish between $Y$ and its  complex conjugate Riemann surface $\bar Y$). In that case, $X$ consists of two connected components.
While it seems artificial at first, the disconnectedness of $G$  is important in our considerations (and this is the reason we work with $\PGLR$ and not with $\SLR$).

The Riemannian metric $g_2$ could be extended to the Riemannian metric $g_3$ on $X$. In fact, it is easier to choose $g_3$ first by choosing  an $Ad(K_0)$-invariant form on the Lie algebra $\mathfrak{g}=sl_2$ (which is isomorphic to the tangent space at any point of $X$), and transporting it to the whole $X$ by the action of $G$.  We also fix the (unique up to a constant) $G$-invariant measure $dx$ on $X$ (e.g., normalized by $\Vol(X)=1$). We denote by $p:X\to Y$ the corresponding projection and by $B_x\simeq S^1\times \bz_2$ the fiber passing through $x\in X$ (in fact, $B_x=xK$ is an $K$-orbit).

Consider a closed geodesic $l\subset Y$ with the natural line element $dl$ on it,  and the $2$-dimensional tori $O=p\inv(l)\subset  X$ ($O\simeq l\times S^1$) obtained as  the preimage of $l$ under $p$ (in general this will consist of two  copies of $2$-tori). Let $T\subset G$ be the subgroup of diagonal matrices. It is well known that $l$ is an orbit of $T$ (in fact, any closed orbit of $T$ gives rise to a closed geodesic on $Y$).
Hence, $O=l\cdot K\subset X$ is the result of the action of the compact subgroup $K$ on the set $l$. Integrating $dl$ with respect to this action, we obtain an area element $do$ on $O$.

The $3$-dimensional manifold $X$ and the hypersurface $O\subset X$ described above are  the geometric ingredients  of our example. In fact, we will  change the setup in the proof in order to avoid certain technical complications explained in Remark \ref{T-invar-funct-sect}. Instead we will deal with another family of tori in $X$. These will also be $S^1$-fibered, but this time over a {\it geodesic circle} on $Y$.  The case of a closed geodesic will be discussed elsewhere. Hence, in this paper, we present the proof for a tori $O=p\inv(\s)\subset X$ fibered over a geodesic circle $\s\subset Y$, $\s=\s_{y_0}(r)=\{y\in Y\ |\ d(y,y_0)=r\}$ centered at an arbitrary point $y_0\in Y$. In fact, we will choose a connected component of this set. The resulting (immersed) hypersurface $O$ clearly is  $K_0$-invariant  (topologically it is  a $2$-dimensional tori).

\subsubsection{Casimir operator} There exists a unique (up to a multiplicative constant) second order differential operator on $X$ which is invariant under the (right) action of $G$ on $X$.  We choose such an operator $C$, called  the Casimir operator on $X$.  In physics it naturally appears in  the Kaluza-Klein theory (see \cite{Z2}).

Let $\partial_b=d/db$ be the unit vector field on $X$ tangential to fibers of the projection $p:X\to Y$. We choose $C$ (or $\Delta$) by requiring that $C=\Delta+\partial_b^2$ (or rather requiring that the $G$-invariant operator $C$ coincides with $\Dl$ on functions constant along fibers of $p$). The symbol of $C$ has the signature $(++-)$. We also will need an elliptic operator on $X$ in order to define Sobolev norms. We will use $\Delta_X=\Delta-\partial_b^2=C-2\partial_b^2$ which coincides (under the appropriate normalization) with the Laplacian on $X$ induced by the Riemannian metric $g_3$.

\subsubsection{Eigenfunctions} The operator $C$ on $X$ has a nice spectral theory (which is admittedly unusual for a hyperbolic operator). Namely, the spectrum of $C$ coincides with that of the elliptic operator $\Delta$ (allowed also to act on forms). However, while all eigenspaces  of $\Delta$ are finite-dimensional, all nontrivial eigenspaces of $C$ are infinite-dimensional.

One  can interpret eigenfunctions of $C$ in terms of the representation theory of $\PGLR$. We denote by $R(g)f(x)=f(xg)$ the action of $G$ on functions on $X$. Let $\Lm\in\br$ be an eigenvalue of $C$. In what follows, it will be more natural (from the point of view of representation theory) to write $\Lm= (1-\lm^2)/4$ for appropriate $\lm\in\bc$. Let $U_\Lm\subset C^\8(X)$ be the space of smooth eigenfunctions of $C$ with the eigenvalue $\Lm$.

\subsection{The restriction}\label{theorem} We fix the set $O\subset X$ as above (i.e., we fix a point $y_0\in Y$, and the {\it radius} $r>0$ of the geodesic circle $\sigma_{y_0}(r)$; $O=p\inv(\sigma_{y_0}(r))$). We are interested in the norm of the natural restriction map on the subspace $U_\Lm$. Namely, consider the restriction  map $r_{O,\Lm}: U_\Lm\to C^\8(O)$, given by  $r_{O,\Lm}(f)=f|_O$ for $f\in U_\Lm$, and
 the corresponding Hermitian form
\begin{equation}\label{H-Lm}H_\Lm(f)=\int_O |r_{O,\Lm}(f)|^2do\end{equation}
defined on the space $U_\Lm$. The form $H_\Lm$ is $K_0$-invariant: $H_\Lm(R(k)f)=H_\Lm(f)$ for all $k\in K_0$, since the set $O$ and the measure $do$ are  $K_0$-invariant.\\

Our main result is
\begin{theorem}\label{claimA}  For any $\eps>0$, there exists a constant $A_{\Lm,\eps}>0$ such that
\begin{equation}\label{3}
 H_\Lm(f)\leq A_{\Lm,\eps}\cdot||f||^2_{W^\eps(X)}\ ,
\end{equation} for any $f\in U_\Lm$.
\end{theorem}
Here we denote by  $W^s(X)$  the $L^2$-Sobolev norm of order $s$ on $X$ (see \cite{BR1} for the representation theoretic treatment).  To define these norms, one can use the elliptic operator $\Delta_X$  on $X$ (e.g., $||f||^2_{W^2(X)}=||f||^2_{L^2(X)}+||\Delta_X(f)||^2_{L^2(X)}$). The value  $s=0$ corresponds to the $L^2$-norm.

Our approach to  Theorem \ref{claimA} is based on the celebrated
 theorem of Gelfand and Fomin \cite{G6}, which interprets  spaces $U_\Lm$  in terms of irreducible unitary representations of $\PGLR$ appearing in  $L^2(X, dx)$ (i.e., each space $U_\Lm$ is a finite direct sum of isomorphic irreducible representations of $G$; see Section \ref{Sect-RepsTh}).  These irreducible representations are called  automorphic representations. This allows one to use powerful methods of the representation theory.
We deduce  the bound \eqref{3} from the corresponding bound for each irreducible component inside of the space $U_\Lm$ (see Theorem \ref{claimB}).

\begin{remarks}\label{rem} 1. The bound \eqref{3} means that on an eigenspace $U_\Lm$ the restriction map $r_{O,\Lm}$ is $L^2-W^\eps$ bounded for any $\eps>0$. We note, however, that the bound we obtain for the constant $A_{\Lm,\eps}$  (i.e., the norm of the restriction map) is quite poor, and corresponds to the general bound for the norm of the restriction of eigenfunctions  to geodesics on $Y$ (i.e., it does not distinguish between positive and negative curvature).
Quantum Chaos intuition suggests the following
\begin{conjecture}\label{conjA} For any fixed $\eps>0$,
\begin{equation*}\label{conjA-bound}
 A_{\Lm,\eps}\ll|\Lm|^\epsilon
\end{equation*} for any $\epsilon>0$.
\end{conjecture}This conjecture is equivalent to the bound \eqref{0} for our particular $X$ and $O$.

2. The space $U_\Lm$ has the natural structure of the tensor product $V_\lm\otimes M_\lm$, where $V_\lm$ is an irreducible (infinite-dimensional for $\Lm\not=0$) unitarizable representation of $G$, and $M_\lm$ is a finite dimensional complex vector space with the dimension $\dim M_\lm=m_\lm$ equal to the {\it multiplicity} of the eigenvalue $\Lm$ of $\Dl$ on $Y$ (one can take $M_\lm$ to be the $\Lm$-eigenspace of $\Delta$ on $Y$). We deal with the restriction norm for a single space $V_\lm$, and them use the fact that $M_\lm$ is finite-dimensional. This is equivalent to  choosing an orthonormal basis in the space $M_\lm$, or what  is the same, choosing an orthonormal basis in the $\Lm$-eigenspace of $\Delta$.

3. A more interesting form to consider is the form which is obtained as the restriction to a single geodesic $l$; namely, the Hermitian form $h_\Lm(f)=\int_l |f|_l|^2dl$ on $U_\Lm$. The form $h_\Lm$ is not $L^2$-bounded. Using a pointwise bound (e.g., as in \cite{BR1}), it is easy to see that it satisfy the bound $h_\Lm(f)\leq A_\eps ||f||^2_{W^{1/2+\eps}(X)}$ for any $\eps>0$ and some $A_\eps$. Forms $H_\Lm$ and $h_\Lm$ are related via integration over $K_0$: $H_\Lm=\int_{K_0} k\cdot h_\Lm dk$, and hence it is not surprising that the form $H_\Lm$ is ``smoother".

4. It seems that our method should provide the following information about the constant $A_{\Lm,\eps}$ (at least for $O$ fibered over a geodesic circle $\s$).

Let $\Delta_B=-\partial_{S^1}^2$ be the Laplacian along $S^1$ fibers  of the bundle $O\to \s$. The operator $\Delta_B$ commutes with $C$ since it corresponds to the action of $\PSOt$ on $X$, and has the simple spectrum $n^2$, $n\in 2\bz$. We denote by $H_\lm$ the form $H_\Lm$ restricted to the (irreducible) automorphic  representation $V_\lm$. Hence for every eigenspace $V_\lm$, we can choose an orthonormal with respect to the form $H_\lm$ basis $\{\phi_\lm^n\}_{n\in 2\bz}$  of joint eigenfunctions of $C$ and $\Delta_B$ ($C\phi_\lm^n=\Lambda \phi_\lm^n$, $\Delta_B\phi_\lm^n=n^2 \phi_\lm^n$). Functions $\phi_\lm^n$ are exponents along the fibers $S^1$. Such a decomposition follows from  representation theory.

To bound the norm of the form $H_\lm$, one needs to bound coefficients $ H_\lm(\phi_\lm^n)$. It is plausible that with more work the following bound could be extracted from our method:
$H_\lm(\phi_\lm^n)\leq  C|n|^\eps,$ for $|n|\geq |\lm|$ and some universal $C>0$.
This would mean that there exists a universal constant in \eqref{3} for $|n|\geq |\lm|$.   For $|n|\leq |\lm|$, we have essentially the standard bound of H\"{o}rmander (or rather the improvement for the restriction to a curved hypersurface) $ H_\lm(\phi_\lm^n)\leq |\lm|^{1/6}$, i.e., $A_{\eps,\lm}\leq |\lm|^{1/6}$ in this range.

One might expect that orthogonal irreducible representations inside the space $U_\Lm$ are ``nearly" orthogonal with respect to the form $H_\Lm$, but we do not know how to prove it. In fact the multiplicities a less problematic for our method than ``clusters" of eigenvalues (and for eigenpackets one expects a different answer).

As a side remark, we note that there is a natural setup of  ``ladders" which allows one to consider eigenfunctions $\phi_\lm^n$ in the framework of PDE theory (in fact, the notion exists on any Riemannian manifold). In particular, Guillemin, Sternberg and Uribe introduced a notion of ``fuzzy ladders'' (see \cite{GS}, \cite{GU},
\cite{ST}, \cite{Z2}). In our setting, this notion concerns a sequence of joint
eigenfunctions of two commuting operators: $\partial_B$ (the Kaluza--Klein operator) and
$\Delta_X$ on $X$. One sees immediately that these joint eigenfunctions are nothing else but
the eigenfunctions $\phi_{\lm}^n$.

We note that there is a variety of effective results  about the asymptotic behavior of eigenfunctions $\phi_\lm^n$ (e.g., effective  equidistribution of $|\phi_\lm^n|^2$ on $SY$ as $|n|\to\8$ shown in \cite{R3}, and uniform bounds for $||\phi_\lm^n||_{L^4(X)}$ for  $|n|\geq |\lm|$ mentioned in \cite{BR3}). These results are obtained through representation theory methods.

5. Results discussed in this paper could be easily extended to noncompact hyperbolic surfaces of finite volume. In fact, this what generated  our original interest in the problem. For such surfaces of the arithmetic origin (e.g., $Y=\operatorname {PSL}_2(\bz)\setminus\CH$), the  Theorem \ref{claimA} is equivalent to a subconvexity bound on certain $L$-functions from Number Theory (see \cite{2}). In particular, the corresponding bound should hold for the Eisenstein series, and also for horocycles.

\end{remarks}

\subsection{Examples related to $S^2$}\label{S2} For some pairs $O\subset X$,  the bound \eqref{0} does not hold. We briefly discuss some examples related to $S^1$ bundles over the standard sphere $S^2$.

First we give a simple example where the local bound \eqref{00} is sharp. Consider $X=S^2\times S^1$, $O=({\rm equator\ of\ }S^2)\times S^1\subset X$,  and $C=\Delta_{S^2}-\frac{1}{4}\Delta_{S^1}$ (i.e., $X$ could be viewed as a wrapped manifold; here $\Delta_{S^n}$ is the (positive) Laplacian on the standard sphere). The hypersurface $O$  is non-characteristic,  time-like and flat. An easy calculation shows that the eigenspace $V_{1/4}$ of $C$ with the eigenvalue $\Lm=1/4$ is infinite-dimensional. An orthonormal basis of $V_{1/4}$ consists of functions $y_{m,l}(s,t)=c\cdot Y^{(m)}_l(s)e^{i(2l+1)t}$, $s\in S^2$ and $t\in S^1$.  Here $Y^{(m)}_l$ are norm one {\it classical} spherical harmonics with the eigenvalue $l(l+1)$ (i.e., $\Delta_{S^2}Y^{(m)}_l=l(l+1)Y^{(m)}_l$, $|m|\leq l$), and $c=1/\sqrt{2\pi}$ is a normalizing constant.  It is easy to see that  the local bound \eqref{00} is achieved in this example for functions in the infinite-dimensional space $V_{1/4}$. This is due to the fact that there are spherical harmonics having maximal possible $L^2$-norm on the equator.

A more relevant example to our discussion  would be the tangent sphere bundle $X= S(S^2)$ over $S^2$. Here the situation is more complicated. Again we  can use  representation theory since $X\simeq \operatorname{SO}(3)$. Hence we choose $X=\operatorname{SO}(3)$ and $G=\operatorname{SO}(3)$ acting on $X$ on the right (we also use the {\it left} action below).   Let  $l\subset S^2$ be a circle on $S^2$ and $O$ its
preimage  in $X$ ($O$ is a $2$-dimensional  tori as before). We can view $l$ as an
orbit  of a  compact subgroup $K'\simeq S^1\subset \operatorname{SO}(3)$ acting on the {\it left} on $X$. Hence $O$ is a double $K'\times \operatorname{SO}(2)$ orbit in $X$ (in such a picture $K'$ could coincide with $\operatorname{SO}(2)$).

There is no {\it natural} hyperbolic operator on X from the point of view of representation theory, but there is a family of such operators which could be described via the right action of $\operatorname{SO}(2)$ on $X$.  Let $C_3$  be the (elliptic) Casimir operator
on $X$  (i.e., the unique (up to a multiple) second  order $\operatorname{SO}(3)$-invariant differential operator on $X$; we normalize it to coincide with  the Laplacian $\Dl_2$ on $S^2$), and let $C_2$ be the differential operator associated with the Casimir operator of $\operatorname{SO}(2)$ (i.e., coming from the  action of $\operatorname{SO}(2)$ on $X$).
We consider a family of operators $C_t=C_3-t\cdot C_2$ depending on the real parameter $t\in \br$. To
obtain  a hyperbolic operator, we have to choose  $t>1$. The distribution of eigenvalues of $C_t$ depends on the size of $t$ and on its  Diophantine property. In particular, for  $t$ which is irrational, one can see that  for big enough $T$, the space $V_T$ which is generated by
eigenfunctions of $C_t$ with the absolute value of the eigenvalue  bounded by $T$ is infinite-dimensional. The same is true for a big enough rational $t$. We can use the  Peter--Weyl theorem on the structure of $L^2(\operatorname{SO}(3))$ in order to describe the space $V_T$. It is not difficult to see that $V_T$ is the direct sum of irreducible representations of $\operatorname{SO}(3)$, and dimensions of these representations are unbounded. The study of the corresponding restriction problem could be reduced to asymptotics of {\it generalized} spherical functions. One can see then that there is an infinite sequence of eigenfunctions in $V_T$ having on $O$ polynomially big (with respect to their  Sobolev norm) restriction norm. However, the exact rate of this growth and its dependence on $O$ and $t$ is more difficult to determine and deserves further investigation.

Finally, we note that the bound \eqref{0} holds for a flat tori $\mathbb{T}^2\subset \mathbb{T}^3$ and  the appropriate (constant coefficients) hyperbolic operator.

 \subsection{Idea of the proof of Theorem \ref{claimA}} The proof of bound \eqref{3} is based on techniques from representation theory. However, the basic idea is quite elementary. We discuss it for a closed geodesic $l$, while in practice we give the proof for a geodesic circle and will discuss closed geodesics elsewhere.

Let $l\subset O\subset X$ be a closed geodesic. In particular, $l=x_0T\subset X$ is a closed orbit for the diagonal  subgroup $T$. Let $V_\lm\subset U_\Lm$ be an irreducible representation of $G$. We denote by $H_\lm$ the form $H_\Lm$ restricted to the automorphic irreducible representation $V_\lm$. We bound the form $H_\lm$, i.e., the form $H_\Lm$ restricted to {\it any} irreducible component inside of $U_\Lm$. The space $U_\Lm$ is  isomorphic, as an abstract representation of $G$, to a finite number of $m_\lm$ copies of $V_\lm$. To obtain the bound for the form $H_\Lm$, we  use the fact that $m_\lm$ is finite.  A simple reduction shows (see section \ref{model-Tfunct}) that in order to bound the norm of the Hermitian form $H_\lm$, it is enough to bound {\it values} of the form $h_\lm(\phi_{e_n})=\int_l |\phi_{e_n}|^2dl$ for some special vectors $\phi_{e_n}\in V_\lm$ (these are vectors which are exponents $e^{in\theta}$ along  fibers $B_x\simeq S^1$; see section \ref{irr-rep}). Here we use the crucial fact that $H_\lm$ (and $H_\Lm$) is $K_0$-invariant, and the structure of irreducible representations of $G$ (i.e., that the space of so-called $K_0$-types are one-dimensional).

 Hence we are interested in values of $h_\lm(\phi_{e_n})=\int_l|\phi_{e_n}|^2dl$ for $|n|\to\8$. We have  $l\simeq S^1$ (in general $l\simeq S^1\cup S^1$ is a union of two closed geodesics, but we disregard this complication since in this paper we will only deal with geodesic circles where this complication could be avoided).   Hence we can use Fourier series expansion on $l$. Consider  expansions $\phi_{e_n}|_l(\te) =\sum_k a_k(\phi_{e_n})e^{ik\te}$ and $h_\lm(\phi_{e_n})=\sum_k|a_k(\phi_{e_n})|^2$.  It turns out that representation theory implies that the coefficients $a_k(\phi_{e_n})$ could be  naturally represented in the form $a_k(\phi_{e_n})=a_k\cdot h(k,n)$, where the function $h(k,n)$ is some special function which is  well approximated by  the classical Airy function.  In particular, from such a representation it follows that for  $|k|\leq (1-\eps)|n|$, we have  $|h(k,n)|^2=O(|n|\inv)$, and that for  $|k|\geq (1+\eps)|n|$, we have $|h(k,n)|^2=O(|n|^{-N})$ for any $N>0$ and any $\eps>0$. In the resonance regime, $k\sim n$, we have the Airy type behavior $|h(k,n)|^2\sim |n|^{-2/3}$. We note that the coefficients $a_k$ are well known in the theory of automorphic functions (see \cite{2}).

 Hence we need to estimate the coefficients $|a_k|^2$.  From simple geometric considerations, one can see that $\sum_{|k|\leq T}|a_k|^2\leq A_\lm T$ for any $T>0$ and some constant $A_\lm>0$ depending on $\lm$. This immediately implies that in the sum $h_\lm(\phi_{e_n})=\sum_k|a_k(\phi_{e_n})|^2=\sum_k|a_k|^2|h(k,n)|^2$, the contribution from regular ranges $|k|\leq (1-\eps)|n|$ and  $|k|\geq (1+\eps)|n|$ is uniformly bounded as $|n|\to\8$.

 Our main point then is that the bound
 \begin{equation}\label{short-sum} \sum_{||k|-T|\leq T^\frac{2}{3}}|a_k|^2\leq B_{\lm,\eps} T^{\frac{2}{3}+\eps} \end{equation}
holds for any $T>0$ and $\eps>0$, and some constant $B_{\lm,\eps} >0$ depending on $\lm$ and $\eps$. Such a  bound clearly implies our main claim \eqref{3} in view of asymptotic for $|h(k,n)|$. To prove \eqref{short-sum}, we use another expansion of $h_\lm$ and the related (less known) representation theory.

 Consider  the following collection of sets: $\Delta l\subset l\times l \subset X\times X$ and $\Delta l\subset \Delta X\subset X\times X$. The Hermitian form $h_\lm$ could be viewed as a linear form on $V_\lm\otimes \bar V_\lm\subset C^\8(X\times X)$ given by the integral over $\Delta l\subset X\times X$.
Clearly, we used the action of $T\times T$ on $ l\times l \subset X\times X$ in order to describe the Fourier decomposition $h_\lm(\phi_{e_n})=\int_l|\phi_{e_n}|^2dl=\sum_k|a_k(\phi_{e_n})|^2$. We can use now the action of $\Delta G$ on $\Delta X\subset X\times X$ in order to give another decomposition of the same form, hence leading to a certain useful identity.

The second decomposition is also easy to describe in elementary terms as follows. Consider the function $|\phi_{e_n}|^2\in C^\8(X)$. We have the spectral decomposition with respect to the  Casimir operator:
$|\phi_{e_n}|^2=\sum _{\lm_i\in Spec(C)}pr_{i}(|\phi_{e_n}|^2)$, where $pr_{i}$ is the projection to the eigenspace $V_{\lm_i}$.
This implies the spectral decomposition $h_\lm(\phi_{e_n})=\int_l |\phi_{e_n}|^2dl=\sum _{\lm_i}\int_l pr_{i}(|\phi_{e_n}|^2)dl$. It turns out that again representation theory allows one to identify quantities $\int_l pr_{i}(|\phi_{e_n}|^2)dl$ in some abstract terms (see Section \ref{triple}), and in particular again allows one to write $\int_l pr_{i}(|\phi_{e_n}|^2)dl=b_{\lm_i}\cdot g(\lm_i,n)$ as a product of an  ``arithmetic" coefficient $b_{\lm_i}$ and a special function $g(\lm_i,n)$.  This shows that we have the identity
\begin{equation}\label{identity}
\sum_k|a_k|^2|h(k,n)|^2 =h_\lm(\phi_{e_n})=\sum _{\lm_i\in Spec(C)}b_{\lm_i}\cdot g(\lm_i,n)\ .
\end{equation}
We use this identity to deduce the bound $\sum_{||k|-T|\leq T^\frac{2}{3}}|a_k|^2\leq B_{\lm,\eps} T^{\frac{2}{3}+\eps}$.

The coefficients $b_{\lm_i}$ satisfy the mean-value bound $\sum_{|\lm_i|\leq T}|b_{\lm_i}|^2\leq CT^2$ (which is in disguise  a bound of L. H\"{o}rmander for the spectral average for a value at a point of eigenfunctions of the Laplacian on $Y$). The function $ g(\lm,n)$ is again some kind of a special function with the Airy type behavior for $|\lm|\sim n$, and typical value of order of $|\lm|\inv$ for non-resonance regime. It seems we are back to the same problem and have gained nothing.

However, and this is the main reason we are able to show some non-trivial saving, the contribution to the Airy behavior of $h(k, n)$ and that of $g(\lm,n)$ comes from {\it different} sets of the ``microlocalization" of the function $\phi_{e_n}$ (the meaning of this in representation theoretic terms is explained in section \ref{automrhic}). In other words, the wavefront sets of distributions defining densities  $h(k,n)$ and $g(\lm,n)$ on the space $V_\lm$ are disjoint. This allows us to construct an appropriate test vector $v\in V_\lm$ such that, on the one hand, $v$ picks up the sum $\sum_{||k|-T|\leq T^\frac{2}{3}}|a_k|^2$ on the left in the identity \eqref{identity}, and, on the other hand, has smooth behavior on the right part of \eqref{identity}.  The mean-value bound on coefficients $b_{\lm_i}$ then implies the bound \eqref{short-sum}.

 \section{Representation theory}\label{Sect-RepsTh}

 \subsection{Gelfand pairs}\label{Gpairs}We will  base our analysis on the notion of  Gelfand pairs (see \cite{GP}). Instead of giving the general definition, we list two cases we need in the proof.

 \begin{itemize}
            \item  For any smooth irreducible representation $V$ of $G=\PGLR$,  the space $\Hom_{K_0}(V,\chi)$ is at most one-dimensional for any character $\chi: K_0\to \bc^\times$ of a maximal connected compact subgroup of $G$ (e.g., $K_0=\PSOt$).
            \item For any three smooth irreducible representations $V_1,\ V_2,\ V_3$ of $G$, the space of trilinear invariant functionals  $\Hom_G(V_1\otimes V_2\otimes V_3, \bc)$ is at most one-dimensional.
          \end{itemize}
In fact, we need only to consider unitarizable  representations. We note that the second multiplicity one statement does not hold for $\SLR$ or $\operatorname{PGL}^+(\br)$ (the space is at most two dimensional). This is the reason we have to deal with the disconnectedness of $\PGLR$.

There is another well-known multiplicity one statement.
For any smooth irreducible representation $V$ of $G$,  the space $\Hom_T(V,\chi)$ is at most one-dimensional for any character $\chi: T\to \bc^\times$ of a maximal tori in $G$ (e.g., for the {\it full} diagonal subgroup in $\PGLR$). This is relevant to the discussion of closed geodesics.

In order to set notations, we quickly review the standard constructions from the theory of automorphic functions and the relevant representation theory (see \cite{G6}, \cite{Bo},\cite{Bu}, \cite{La}).

\subsection{Irreducible representations of $\PGLR$}\label{irr-rep} Infinite-dimensional  irreducible unitary representations of $G$ are naturally split into two types: induced representations and discrete series. Unitary induced representations $\pi_{\lm,\eps}$ are parameterized  by a pair $(\lm,\varepsilon)$, where $\lm$ is a complex number which belongs to the set $\in i\br\cup(0,1)$ and $\varepsilon\in\{0,1\}$. Discrete series representations $\pi_k$ are parameterized by an even positive integer $k\in 2\bz^+$. The parameter $\Lm$ we used before describes the action of a particular generator of the center of the universal enveloping algebra $U(sl_2)$ (i.e., the Casimir operator). The relation between two parameters  is $\Lm=(1-\lm^2)/4$ for induced representations, and $\Lm=k/2(1-k/2)$ for the integer case. Induced representations parameterized by $i\br$ are called principal series representations, and  those parameterized by $(0,1)$ are called complementary series representations. We call these even for $\eps=0$, and odd for $\eps=1$ (although this deviates from the standard terminology where this term is used in connection with the central character of $\GLR$). Unitary induced representations and discrete series representations, together with two one-dimensional representations, $\det^\eps$, $\varepsilon\in\{0,1\}$, exhaust all irreducible unitary representations of $G$.

We note that only representations of principal and complementary series representations have a non-zero $K_0$-fixed vector, and hence give rise to eigenfunctions of $\Delta$ on $Y$.  Representations of the complementary series correspond to what is called the exceptional spectrum of $Y$ (i.e. $\Lm<1/4$). Representations of the discrete series correspond to forms on $Y$. We will assume for simplicity that there is no exceptional spectrum, i.e., that all automorphic representations with non-zero $K_0$-fixed vectors are representations of the unitary principal series.

The structure of  unitary representation of $\PGLR$ is well known. Induced  representations could be modeled in various spaces of functions on various manifolds. In particular, we have the following realization of the space $V_{\lm,\eps}$ of smooth vectors for the induced representation $\pi_{\lm,\eps}$. Consider the space  $\CH_{\lm,\eps}$  of smooth even homogeneous functions on $\br^2\setminus 0$ of the homogeneous degree $\lm-1$ (i.e., $f(tv)=|t|^{\lm-1}f(v)$ for any $t\in \br^\times$ and $0\not=v\in\br^2$). The evenness  condition is necessary since we consider the group $\PGLR$. We have the natural action of $\GLR$ given by $\pi_{\lm,\eps}(g)f(v)=f(g\inv v)\cdot |\det(g)|^{(\lm-1)/2}\cdot \det(g)^\eps$, which is trivial on the center and hence defines a representation of $\PGLR$. We call such a realization the plane realization. We will use the realization of such a representation in the space $C_{ev}^\8(S^1)$  of  smooth {\it even} functions on $S^1$ (i.e., $f(\theta+\pi)=f(\theta)$) which is obtained by taking the restriction to the circle $S^1\subset \br^2\setminus 0$.  We call such a realization  the circle model. In such a realization, the action of $K_0\simeq S^1$ is given by the rotation. The invariant unitary norm then coincides with the standard norm on $L^2(S^1)$ (this is where the assumption $\lm\in i\br$ is used; for complimentary series the norm is not local). Hence there is a natural orthonormal basis $\{e_n(\te)=e^{in\te}\}_{n\in2\bz}$ consisting of $K_0$-equivariant vectors (called $K_0$-types). The parameter $\eps$ specifies the action of the element $\dl= \left(
                                          \begin{smallmatrix}
-1 &  \\
                                            & 1 \\
\end{smallmatrix}
                                        \right)$ with the negative determinant. Namely, the action of $\dl$ is given by $\pi_{\lm,\eps}(\dl)f(\theta)=(-1)^\eps f(-\theta)$. In particular, $\pi_{\lm,\eps}(\dl)e_n=(-1)^\eps e_{-n}$ for the basis of $K_0$-types.

The discrete series representations lack such a simple geometric model (a source of many computational difficulties) , and we discuss necessary amendments in \ref{DiscrSer}.

\subsection{Automorphic representations}\label{automrhic} The notion of automorphic representation allows one to use effectively models of representations in our setup. Namely, till now we viewed the space $V_\lm\subset C^\8(X)$ as a subspace of functions on $X$. This is a rather inaccessible (however remarkable) realization of an irreducible representation. For principal series, the alluded Gelfand and Fomin theorem  implies that there is a $G$-equivariant map $\nu_\lm: C_{ev}^\8(S^1)\to V_\lm$, which we can assume to be an isometry. This means that every vector $\phi\in V_\lm$ is of the form $\phi_u=\nu_\lm(u)$ for some function $u\in C_{ev}^\8(S^1)$. This allows us (sometimes) to translate various questions about eigenfunctions of $C$ on $X$ (or of $\Delta$ on $Y$) into questions about $\nu_\lm$ and an appropriate vector $u\in C_{ev}^\8(S^1)$. We may argue that the correspondence $\phi\leftrightarrow u$ is somewhat analogous to what is called microlocalization of eigenfunctions.

For discrete series, similar considerations apply. However, discrete series representations lack a good realization similar to the circle model above. We will discuss the appropriate changes in \ref{DiscrSer}.

We assumed that the space $X$ is compact (i..e, the Riemann surface $Y$ is compact). Let
\begin{eqnarray}\label{L2X-spec}L^2(X)=\Bigl(\oplus_i(L_{\tau_i,\eps_i},\nu_{i})\Bigr)\oplus\Bigl(\oplus_j
(L_{k_j},\nu_{j})\Bigr)\end{eqnarray} be the decomposition into
the orthogonal sum of irreducible unitary representations of $G$. Here $\nu_i:L_i\simeq
L_{\tau_i,\eps_i}\to L^2(X)$ are unitary representations of class one
(i.e., those which correspond to Maass forms on $Y$ with the eigenvalue $(1-\tau_i^2)/4$), and
$\nu_{j}:L_{k_j}\to L^2(X)$ are representations of discrete series (i.e., those
which correspond to holomorphic forms on $Y$). We denote by
$V_i\subset L_i$ the corresponding spaces of smooth vectors and by
${\rm pr}_{L_i}:L^2(X)\to \nu_i(L_i)$ the corresponding orthogonal
projections (note that ${\rm pr}_{L_i}:C^\8(X)\to V_i$). Note that there might be multiplicities for automorphic representations (i.e., two realizations $\nu,\ \nu':L_\pi\to  L^2(X)$ of an irreducible representation $(L_\pi,\pi)$, having different images) notably for representations of discrete series.

\subsubsection{$\PGLR$ versus  $\PSLR$}\label{PGL-PSL} A more familiar setup is that of automorphic functions on $\PSLR$. We describe how it is connected to our setup. Let $G_1=\PSLR$ and $\G_1\subset G_1$ be a (co-compact) lattice. We consider the quotient space $X_1=\G_1\setminus G_1\simeq SY$ which is isomorphic to the tangent sphere bundle over a (compact) Riemann surface $Y$. The spectral decomposition of the Casimir operator on $X_1$ (or of the Laplacian $\Dl$ on $Y$) leads to the decomposition $L^2(X_1)=\oplus_i(N_i,\s_i)$ into irreducible unitary representations of $\PSLR$. Now consider $\G_1$ as a lattice in $G=\PGLR$ (or in $\operatorname{PSL}^\pm_2(\br)$). Formally,  the space $X=\G_1\setminus G=X_+\cup X_-$ consists of two connected components. It might happen however that there exists a lattice $\G\subset G$ such that $\G_1\subset \G$  and $\G\cap \operatorname{PGL}_2^-(\br)\not=\emptyset$. The best known example is $\G_1=\operatorname{PSL}_2(\bz)$ and $\G=\operatorname{PGL}_2(\bz)$. A geometric example could be constructed from a bordered Riemann surfaces $Y$ with a totally geodesic boundary, by taking the  connected sum of $Y$ and $\bar Y$ along the geodesic boundary.

We want to connect representations of $G$ in $L^2(X)$ and those of $G_1$ in $L^2(X_1)$.

First let us consider the situation where $\G=\G_1$. Let $\dl= \left(
                                          \begin{smallmatrix}
-1 &  \\
                                            & 1 \\
\end{smallmatrix}
                                        \right)$.
We have $X=X_+\cup X_-$ with $X_+=\G_1\setminus G_1$, and $X_-\simeq X_1$ under the multiplication by $\dl$, i.e.,  $x\mapsto x\dl$.

Let $(N,\s)\subset L^2(X_1)$ be an automorphic  representation of $G_1$  which as an abstract representation is induced (i.e., of principal or complimentary series). Consider the space $N'\subset L^2(X_-)$ of functions obtained by the action of $\dl$ (i.e., $N'=\{f'\ | \ f'(x\dl)=f(x)\ {\rm for\ some\ }\ f\in N\}$). We have then $N\oplus N'=L_0\oplus L_1$, where $L_\eps$, $\eps=\{0,1\}$ are (even and odd) irreducible representations of $G$ which as representations of $G_1$ are isomorphic. Under the restriction map $L^2(X)\to L^2(X_1)$, these give the same space of eigenfunctions of the Casimir on $X_1$.

For an automorphic representation $(N, \s)$ of $G_1$ of discrete series,  one first considers the complex conjugate representation  $(\bar N, \bar \s)$ (i.e., the space of functions $\bar f$, $f\in N$ with the natural action of $G_1$), and then the action of $\dl$ (interchanging connected components). In that case, we obtain two copies of isomorphic representations of $G$, and each of these split into two irreducible components under the restriction to $G_1$ (holomorphic and anti-holomorphic discrete series representations on different components of $X$).

Taking restrictions  of functions on $X$ to subsets $X_\pm\simeq X_1$, we obtain two isomorphisms $i_\pm:L^2(X_\pm)\simeq L^2(X_1)$ as representations of $G_1$. Note  that $i_\pm$ are {\it not} coming from the representation-theoretic restriction of representations of $G$ to representations of $G_1$.

The case of $X$ consisting of one connected component is more interesting.  Geometrically this means that the {\it oriented} hyperbolic Riemann surface $Y$ is isometric (as an oriented Riemann surface) to itself after reversing the orientation. We have $X_1=X$, and $\dl$ acts on it on the right.  Let us assume for simplicity that $\dl$ normalizes $\G_1$. Consider an automorphic irreducible unitary representation $(L,\pi)\subset L^2(X)$ of $G$. If $\pi=\pi_{\lm,\eps}$ is induced, then we can take its restriction to the representation of $G_1$ which is again irreducible. This leads to even and odd eigenfunctions (with respect to $\dl$) of $C$ on $X$ depending on the parity of $\eps$. Note that we can not reverse the direction by starting with some automorphic  representation of $G_1$ and extending action of $G$ to it  (if there are multiplicities, one has to split the space into the sum of $\dl$-even and $\dl$-odd subspaces).

For an automorphic discrete series representation of $G$, the restriction to $G_1$ leads to two irreducible automorphic components which are interchanged by $\dl$.

\subsection{Geodesic circles}\label{K-geometry}We start with
the geometric origin of the spherical Fourier coefficients.
We fix the ``standard"  maximal compact subgroup $K=\POt\subset G$, and the
identification $G/K\to\CH$, $g\mapsto g\cdot i$, with the hyperbolic plane. Let $\G\subset \PSLR$ be as before. We denote by $p_\CH:\CH\to \G\sm \CH\simeq Y$ the corresponding
projection. It is compatible with the distance function $\rm d(\cdot,\cdot)$ on $Y$
and $\CH$. Let $y\in Y$ be a
point, and let $R_y>0$ be the injectivity radius of $Y$ at $y$.
For any $r< R_y$, we define the geodesic circle of radius $r$
centered at $y$ to be the set $\s(r,y)=\{y'\in Y|{\rm
d}(y',y)=r\}$. Since the map $p_\CH$ is a local isometry, we have that
$p_\CH(\s_\CH(r,z))=\s(r,y)$ for any $z\in\CH$ such that $p_\CH(z)=y$,
where $\s_\CH(r,z)$ is the corresponding geodesic circle in $\CH$
(all geodesic circles in $\CH$ are the Euclidian circles in $\bc$, though
with a different center from $z$). We associate to any such circle
on $Y$ an orbit of a compact subgroup on $X$. Namely, let
$K_0=\operatorname{PSO(2)}\subset K$ be the connected component of $K$. Any
geodesic circle on $\CH$ is of the form $\s_\CH(r,z)=hK_0g\cdot i$
with $h,\ g\in G$ such that $h\cdot i=z$ and $hg\cdot i\in
\s_\CH(r,z)$ (i.e., an $h$-translation of a standard geodesic
circle centered at $i\in\CH$ and passing through $g\cdot
i\in\CH$). Note that the radius of the circle is given by the
distance ${\rm d}(i,g\cdot i)$ and hence $g\not\in K_0$ for a
nontrivial circle. Given the geodesic circle $\s(r,y)\subset Y$, we
consider a circle $\s_\CH(r,z)\subset \CH$ projecting onto
$\s(r,y)$ and the corresponding elements $g,\ h\in G$. We denote
by  $K_\s=g\inv K_0g$ the corresponding compact subgroup and
consider its orbit $\ck_\s=hg\cdot K_\s\subset X$. Clearly we have
$p(\ck_\s)=\s$. We endow the orbit $\ck_\s$ with the unique
$K_\s$-invariant measure $d\mu_{\ck_\s}$ of the {\it total mass
one} (from the geometric point of view a more natural measure
would be the length of $\s$).

We note that in what follows the restriction $r<R_y$ is not
essential. From now on, we assume that $\ck\subset X$ is an orbit
of a {\it maximal} connected compact subgroup $K'_0\subset G$. The group $K'_0$ is conjugate to $\PSOt$, but  does not coincide with it. We denote by $K'$ the corresponding maximal compact subgroup in $G$. Let $\ck\subset  X$ be an orbit of $K'_0$.  Since we have assumed that $X$ consists of two components $X=X_+\cup X_-$, we also have $\ck\cdot K'=\ck_+\cup \ck_-$, $\ck_\pm\subset X_\pm$ and each of these sets is an orbit of $K'_0$. We will assume that $\ck=\ck_+\subset X_+$. The restriction $r<R_y$ simply means that
the projection $p(\ck)\subset Y$ is a smooth non-self intersecting
curve on $Y$. We remark that  polar
geodesic coordinates $(r,\theta)$ centered at a point
$z_0\in\CH$ could be obtained from the Cartan
$KAK$-decomposition of $\PSLR$ (see \cite{He} where  analysis on $\CH$ is discussed in the context of the representation theory).

\subsection{ Generalized periods and equivariant
functionals}\label{K-fourier}  We fix a point $\dot{o}\in\ck$. To
a character $\chi:K'_0\to \bc^\times$, we associate a function
$\chi_.(\dot{o}k')=\chi(k')$, $k'\in K'_0$ on the orbit $\ck$, and
the corresponding functional on $C^\8(X)$ given by the generalized period
\begin{align}\label{O-funct-K}
d^{aut}_{\chi,\ck}(f)=\int_\ck f(k)\bar{\chi}_.(k)d\mu_\ck
\end{align}
for any $f\in C^\8(X)$. The functional $d^{aut}_{\chi,\ck}$ is
$\chi$-equivariant with respect to the right action of $K'_0$ on functions on $X$:
$d^{aut}_{\chi,\ck}(R(k')f)=\chi(k')d^{aut}_{\chi,\ck}(f)$ for any
$k'\in K'_0$ (here $R$ is the right action of $G$ on the space of
functions on $X$). For a given orbit $\ck$ and a choice of a
generator $\chi_1$ of the cyclic group $\hat{K'_0}\simeq \bz$ of characters of
the compact group $K'_0$, we will use the shorthand notation
$d^{aut}_{n}=d^{aut}_{\chi_n,\ck}$, where $\chi_n=\chi_1^n$. The
functions $(\chi_n)_.$ form an orthonormal basis for the space
$L^2(\ck,d\mu_\ck)$.

 Let $\nu:V\to C^\8(X)$ be an irreducible automorphic
representation. When it does not lead to confusion, we denote by
the same letter the functional
$d^{aut}_{\chi,\ck}=d^{aut}_{\chi,\ck,\nu}$ on the space $V$
induced by the functional $d^{aut}_{\chi,\ck}$ defined above on
the space $C^\8(X)$. Hence we obtain an element in the space
$\Hom_{K'_0}(V,\chi)$. We next use the above mentioned multiplicity one property, i.e., the fact that
the pair $(G,K'_0)$ is a Gelfand pair (i.e., that $\dim \Hom_{K'_0}(V,\chi)\leq 1$).

Let $V\simeq V_{\lm,\eps}$ be a representation of the principal series.
We have  then the multiplicity one property $\dim\Hom_{K'_0}(V_{\lm,\eps},\chi)\leq 1$ for any character
$\chi$ of $K'_0$ (i.e., the space of $K'_0$-types is at most one
dimensional for a maximal connected compact subgroup of $G$). In
fact, $\dim\Hom_{K'_0}(V_{\lm,\eps},\chi_n)=1$ if and only if $n$ is even.

To construct a model $\chi$-equivariant functional on $V_{\lm,\eps}$, we
consider the circle model $V_{\lm,\eps}\simeq C^\8_{ev}(S^1)$ in the
space of even functions on $S^1$ and the standard vectors
(exponents) $e_n=\exp( in\theta)\in C^\8(S^1)$ which form a basis
of $K_0$-types for the {\it standard} connected maximal compact subgroup
$K_0=\PSOt$. For any $n$ such that
$\dim\Hom_{K_0}(V_{\lm,\eps},\chi_n)=1$, the vector
$e'_n=\pi_\lm(g\inv)e_n$ defines a non-zero
$(\chi_n,K'_0)$-equivariant functional on $V_{\lm,\eps}$ by the formula
\begin{eqnarray}\label{d-mod-matrix-c}
d^{mod}_{n}(v)=d^{mod}_{\chi_n,\lm}(v)=\langle v,e'_n\rangle\ .
\end{eqnarray} We call such a functional the {\it model} $\chi_n$-equivariant
functional on the representation $V\simeq V_{\lm,\eps}$.

The uniqueness principle implies that there exists a constant
$a_n=a_{\chi_n,\ck}(\nu)\in \bc$ such that
\begin{eqnarray}\label{b-def}
d^{aut}_{n}(u)=a_n\cdot d^{mod}_{n}(u)\ ,
\end{eqnarray} for any $u\in V$. Here we suppressed the dependence on the orbit and the automorphic representation since these are fixed in our discussion (nevertheless this dependence is central in other applications).

\subsubsection{Invariant functional}\label{T-invar-funct-sect}  Since the $K'_0$-invariant functional will play a double role in our construction, we introduce another notation for it. Note that such a functional automatically vanishes on a representation of discrete series.

Let $\chi_0\equiv 1$ be the trivial character of $K'_0$.  We denote by $d_{\tau,\eps}(u)=\langle u,e'_0\rangle_{V_{\tau,\eps}}$, $u\in  V_{\tau,\eps}$, the corresponding  model
functional. Let  $\nu_i:V_{\tau_i,\eps_i}\to C^\8(X)$ be an irreducible automorphic representation
of class one. We have as before
\begin{eqnarray}\label{T-invariant-f}
d_i^{aut}(u)=d^{aut}_{\chi_0,\ck,\nu_i}(u)=\int_\ck
\nu_i(u)(k)\bar{\chi_0}_.(k)d\mu_\ck=\alpha(i)d_{\tau_i,\eps_i}(u)\
,
\end{eqnarray}
for any $u\in V_{\tau_i,\eps_i}$, and a constant $\alpha(i)=\alpha_\ck(\nu_i)\in\bc$.

We want to compare the coefficients $\alpha(i)$ with more
familiar quantities.  Let $\ck=x_0\cdot K'_0\subset X$ and let  $\phi'_{\tau_i}=\nu_i(e'_0)$ be the
automorphic  function which corresponds to a $K'_0$-invariant vector
$e'_0\in V_{\tau_i,\eps_i}$ of norm one. From the definition of
$d_{\tau,\eps}$, it follows that
\begin{eqnarray}\label{b-0-nomalize}
\alpha(i)=\phi'_{\tau_i}(x_0)\ .
\end{eqnarray}
Hence, under the normalization that we choose, the coefficients
$\alpha(i)$ are equal to values at a point $x_0$ for Maass forms on
the Riemann surface $Y'=\G\sm G/g\inv Kg$.

\begin{remark}\label{abcent-discr-ser} A trivial, but important, remark is that on the discrete series representations any
$K'_0$-invariant functional is identically zero, unlike a $T$-invariant functional for the diagonal subgroup $T$ of $G$. This greatly simplifies the technicalities in what follows; and in fact, this is the reason that we treat only {\it geodesic circles} and leave {\it closed geodesics} aside. For the later, we need to discuss discrete series representations at length, and the resulting computations are more involved.
\end{remark}

\subsection{First Gelfand pair: abelian spectral decomposition on $\s$}\label{model-Tfunct} Let $\nu: V_{\lm,\eps}\to L^2(X)$ be an automorphic representation of principal series.
Consider the Hermitian form on $V_{\lm,\eps}$ arising from the restriction to the orbit $\ck$ of $K'_0$  (i.e., to a geodesic circle $\s$) \begin{equation}h_\lm(\phi_u)=\int_\ck\left|\phi_u\bigm|_\ck\right|^2d\ck, \end{equation}
for $\phi_u=\nu(u)$, $u\in V_{\lm,\eps}$. Obviously the Hermitian form $h_\lm$ could be decomposed with respect to the action of $K'_0$.  We have
\begin{equation}\label{spec-h-lm}
h_\lm(\phi_u)=\sum_k|d^{aut}_k(\phi_u)|^2=\sum_k|a_k|^2|d^{mod}_k(u)|^2\ .
\end{equation} This decomposition could be used in order to evaluate the integrated form \begin{equation}H_\lm=\int_{K _0}k\cdot h_\lm\ dk\ .\end{equation}
This is the Hermitian form which appears in Theorem \ref{claimA}, restricted to an irreducible component inside of the space $U_\Lm$.

The form $H_\lm$ is $K_0$-invariant, and hence is determined by its values on the (essentially unique) orthogonal (with respect to $H_\lm$) basis of $V_{\lm,\eps}$ consisting of $K_0$-equivariant vectors $\{e_n\}_{n\in 2\bz}$  described in section \ref{irr-rep}. Hence we need to estimate the quantities $H_\lm(\phi_{e_n})=h_\lm(\phi_{e_n})$, where we denote  $\phi_{e_n}=\nu_\lm(e_n)$. Substituting  \eqref{spec-h-lm} into the expansion \eqref{b-0-nomalize}, we see that
\begin{equation}\label{specH-lm-aut}
H_\lm(\phi_{e_n})=h_\lm(\phi_{e_n})=\sum_{k\in 2\bz}|a_k|^2|d_k^{mod}(e_n)|^2\ ,
\end{equation} for any $n\in 2\bz$.

In order to prove bound \eqref{3}, we need to show that $H_\lm(\phi_{e_n})\ll |n|^\dl$ for any $\dl>0$ and fixed $\lm$.

\begin{theorem}\label{claimB}  For any $\dl>0$, there exists a constant $A_{\lm,\dl}>0$ such that
\begin{equation}\label{4}
 \sum_k|a_k|^2|d_k^{mod}(e_n)|^2\leq A_{\lm,\dl} |n|^\dl\ ,
\end{equation} for any $n\in 2\bz$.
\end{theorem}

This implies bound \eqref{3} and proves the main theorem of the paper. Below we present the proof for the principal series representations, and discuss in \ref{DiscrSer} amendments needed to be made for the discrete series.

\begin{remarks}{}  1. The coefficients $a_k$ are well known in the theory of automorphic functions. These are called the spherical Fourier coefficients of Maass forms and were introduced by H. Peterson. Quantities  $|a_k|^2$ are related to special values of $L$-functions (for special geodesic circles on Riemann surfaces of arithmetic origin, and for what is called Hecke-Maass forms), and are of  utmost importance in Number Theory.  In particular, it is generally believed that they satisfy the bound $|a_k|\ll |k|^\eps$ for any $\eps>0$ (the Lindel\"{o}ff conjecture).  We stress that coefficients $a_k$ depend on $\s$, $\nu$ and  $\G$.

2. It is relatively easy to prove the mean-value bound
\begin{equation}\label{mean}\sum_{|k|\leq T}|a_k|^2\leq AT\ ,
 \end{equation}
 for any $T\geq |\lm|$, and some universal constant $A$ depending on $Y$ only (this follows from \cite{BR2}, and was spelled out explicitly in  \cite{1}).\end{remarks}

\subsection{Oscillating integrals} We treat principal series representations first and discuss discrete series in \ref{DiscrSer}. The coefficients  $d_k^{mod}(e_n)$  should be viewed as a function of two variables $k$ and $n$ (in fact, there is also a dependence on the parameter of representation $\lm$ which we suppress as  it is fixed in our discussion). Recall that we defined these by the matrix coefficient \eqref{d-mod-matrix-c}
\begin{eqnarray}\label{d-mod-matrix-c2}
d^{mod}_{k}(e_n)=\langle e_n,\pi_\lm(g\inv)e_k\rangle\ .
\end{eqnarray} Here  $e_j\in V_{\lm,\eps}$ is the norm one $j$-$K_0$-type given by the function $e_j(\theta)=e^{ij\theta}$ in the circle model $V_{\lm,\eps}\simeq C^\8_{even}(S^1)$.
Such a matrix coefficient clearly could be computed via an  oscillating integral in one variable. Taking into account the action of $K_0$ in the space $V_{\lm,\eps}$, we have
\begin{align}\label{osc-int}
 d_k^{mod}(e_n)=\int_{S^1}|g'(\theta)|^{-\haf}
e^{\haf\lm\ln|g'(\theta)|}e^{ikg(\theta)-in\theta}d\theta\ .
\end{align}

Note that this formula does nor depend on the parity $\eps$ of $V_{\lm,\eps}$. This is not surprising as $\eps$ describes the action of the element $\dl$.

Integrals \eqref{osc-int} are common in the theory of special functions, and  have been  well-studied in  Classical Analysis \`{a} la Whittaker and Watson. In fact, these are well known Legendre functions.
The asymptotic of such an integral is controlled by the  stationary phase method.  We are interested in the behavior of   $d_k^{mod}(e_n)$ for large parameters $k$ and $n$, for  $\lm$ and $g$ which are fixed. Hence we denote by $A(\theta)=|g'(\theta)|^{-\haf}e^{\haf\lm\ln|g'(\theta)|}$ the amplitude and by $S_{k,n}(\theta)=kg(\theta)-n\theta$ the phase in the oscillating integral \eqref{osc-int}.   By computing the relevant phase function, one can easily see that the function $d_k^{mod}(e_n)$ is well modeled by the classical Airy function $\mathbb{A}(x)$ (see  \cite[Theorem 7.7.18]{He}; we note that careful analysis of such a reduction  is done in \cite{BR3}).

Let $M=\max_{\theta\in S^1}|g'(\theta)|>1$. We have $M\inv\leq g'(\theta)\leq M$ as $g$ is a linear-fractional map on  $S^1$.
We have three types of behavior for the integral  \eqref{osc-int}.

\begin{enumerate}
  \item For $k$ in the range $1.1M\inv n\leq k\leq 0.9 M n$, the phase function  has two critical points of Morse type;
in this case we can estimate the integral using the stationary phase
method. In this range the integral is of order of $|k|^{-\haf}$.
 \item For $k$ close to $M\inv n$ or to $M n$, critical points of the phase collide to a cubic critical point. In order to get uniform bounds in this region, we can use properties of the Airy function.
\item For $k$ in ranges $k\leq 0.9 M\inv n$ or $k\geq 1.1 M n$,
there are no critical points of the phase, and the value of the integral is negligible (in $n$). This region could be discarded in further analysis.

\end{enumerate}

We treat the case  $k\asymp Mn$, and the complimentary case $k\asymp M\inv n$ is identical. Namely, we assume that $n>0$ and $k\geq n$ (the treatment for the range $k<n$ is identical).  Let us introduce the notation $\dl=\dl(k,n)=Mn/k$.

The above analysis (made rigorously in \cite{2} and \cite{BR3}) implies  in particular that for any fixed sufficiently small $\eps>0$, in the region  $1+\eps\geq \dl \geq1-\eps $  (i.e., $k\asymp M n$), we have
\begin{eqnarray}\label{bddd}|d^{mod}_{k}(e_n)-A(k,\dl) | \leq  C   |k| ^{-\frac{2}{3}}\ . \end{eqnarray}
Here $A(k,\dl) = |k|^{-\frac{1}{3}}\mathbb{A}(k^{\frac{2}{3}}(\dl-1))=|k|^{-\frac{1}{3}}\mathbb{A}(k^{-\frac{1}{3}}(Mn -k))$, where $\mathbb A$ is  the classical Airy function (see  \cite{He}). We only need to know that the Airy function $\mathbb{A}(x)$ is a smooth function, bounded at $0$, rapidly decaying for $x>0$ and uniformly bounded by $|x|^{-1/4}$ at infinity.
This well-known asymptotic of the Airy function imply that (roughly) we have three types of behavior:
\begin{description}
  \item[\it Regular:] If  $M-1\leq (Mn-k)/k\leq 1-\eps$, then $|d_k^{mod}(e_n)|\leq |n|^{-\haf}$.
  \item[\it Resonance:] If $1-\eps\leq (Mn-k)/k\leq 1+\eps$, then $ |d_k^{mod}(e_n)|\leq |k|^{-\frac{1}{3}}(1+|k|^{2/3}|Mn-k|)^{-1/4}$.
  \item[\it Cutoff:] If $(Mn-k)/k\geq1+\eps$, then  $|d_k^{mod}(e_n)|\ll |k|^{-N}$ for any $N>0$.
\end{description}

Splitting the summation in \eqref{4} according to these regions of $k$, we see that
the mean-value bound \eqref{mean} allows one to treat the regular part of the summation in \eqref{4}, i.e., $k$ such that $\dl$ is not close to $1$. This region  of summation in $k$ is of  length comparable to $n$ due to the cutoff, and the value of $|d_k^{mod}(e_n)|^2$ is of order of $|n|\inv$. Hence we see that the total contribution from this region of $k$ is {\it uniformly} bounded as $n\to\8$.

We are left with the ``short" sums of the form $\sum_{|\dl(k,n)-1|\ll T\inv} |a_k|^2$ for $n$ and $T\to\8$, $T\ll n$. Denoting $T=Mn$, we consider  a  sum  $\sum_{||k|-T|\leq |T|^\gamma}|a_k|^2$ for $\g<1$. The size of the transition region of the Airy function implies that we only need to consider  the range $1\geq\g\geq2/3$ (in fact, we do not know how to deal with shorter sums!).
It turns out that these sums could be  bounded effectively with the help of the second Gelfand pair  (i.e., using  the multiplicity one  statement $\dim \Hom_G(V_1\otimes V_2\otimes V_3, \bc)\leq 1$). This is the more ``tricky" part of the proof and is done in \cite{2}. We explain it in the next section.

\section{Triple product}\label{triple}
 Here we explain how to obtain the bound
\begin{eqnarray}\label{short-bound-an}\sum_{||k|-T|\leq T^\frac{2}{3}}|a_k|^2\leq B_{\lm,\eps} T^{\frac{2}{3}+\eps}\ . \end{eqnarray} This is the content of Theorem 1.5 from \cite{2}. We recall the corresponding setup.
\subsection{Second Gelfand pair: triple product spectral decomposition} Recall that we consider the Hermitian form given by $h_\lm(\phi)=\int_\ck|\phi_u|_\ck|^2d\ck$ for $u\in V_{\lm,\eps}$. We now switch to the corresponding linear functional on $E=V_{\lm,\eps}\otimes  V_{-\lm,\eps}$ which we denote by the same letter. We have $h_\lm(w)=\int_\ck r_\Dl(\nu_E(w))d\ck$ for $w\in E$. We now  consider the spectral expansion for the form $h_\lm$ coming from triple products.

\subsubsection{Spectral theory on $X$} Let $\nu:V\to C^\8(X)$ be an irreducible automorphic
representation as before and $\nu_E=\nu\otimes\bar\nu: E=V\otimes
\bar V\to C^\8(X\times X)$ the corresponding realization of $E$.
We assume that the space $X$ is compact, and hence we have the discrete sum decomposition \eqref{L2X-spec}
\begin{eqnarray}\label{L2X-spec2}L^2(X)=\left(\oplus_i (L_i,\nu_i)\right)\oplus\left(\oplus_\kappa
(L_\kappa,\nu_\kappa)\right)\end{eqnarray} into
irreducible unitary representations of $G$. Here $\nu_i:L_i\to L^2(X)$ are unitary representations of class one
(i.e., those which correspond to Maass forms on $Y$ with the eigenvalue $(1-\tau_i^2)/4$) and
$L_\kappa$ are representations of discrete series (i.e., those
which correspond to holomorphic forms on $Y$). We denote by
$V_i\subset L_i$ the corresponding spaces of smooth vectors and by
${\rm pr}_{L_i}:L^2(X)\to \nu_i(L_i)$ the corresponding orthogonal
projections (note that ${\rm pr}_{L_i}:C^\8(X)\to V_i$).

We consider triple products of eigenfunctions. Let
$r_\Dl:C^\8(X\times X)\to C^\8(X)$ be the map induced by the
imbedding $\Dl: X\to X\times X$. Let $\nu_i:V_{\tau_i,\eps_i}\to C^\8(X)$
be an irreducible automorphic representation. Composing $r_\Dl$
with the projection ${\rm pr}_{L_i}:C^\8(X)\to\nu_i(V_{\tau_i,\eps_i})$, we
obtain  $\Dl G$-invariant map $T^{aut}_i:E\to
V_{\tau_i,\eps_i}$ and the corresponding automorphic trilinear functional
$l^{aut}_{i}$ on $E\otimes V_{\tau_i,\eps_i}^*$ defined by
$l^{aut}_{i}(v\otimes u\otimes t)=\langle r_\Dl(\nu_E(u\otimes v)),\bar
\nu_i(t)\rangle$ (here we identified $\bar V_{\tau_i,\eps_i}$ with the smooth
part of $V^*_{\tau_i,\eps}\simeq V_{-\tau_i,\eps_i}$). Such a functional is clearly
$G$-invariant, and hence we can invoke the uniqueness principle
for trilinear functionals (see section \ref{Gpairs}).

\subsubsection{Triple product spectral expansion}
Recall that  in section \ref{T-invar-funct-sect} we denoted by  $d^{aut}_i:V_i\to \bc$ the $K'_0$-invariant functional coming from the integration along the orbit $\ck\subset X$ (see \eqref{T-invariant-f}). The expansion \eqref{L2X-spec2} implies that
\begin{eqnarray}\label{h-lm-spec} h_\lm(w)=\sum_i d^{aut}_i(T^{aut}_i(w))\ ,
\end{eqnarray} for any $w\in E$.
 We now write $d^{aut}_i(T^{aut}_i(w))$ as an integral transform in the circle model of representations.

\subsubsection{Model triple product} Let $V_{\lm,\eps}$ and $V_{\tau,\eps'}$ be two irreducible unitary representations of principal series (and we denote by $E$ the smooth part of the tensor product  $V_{\lm,\eps}\otimes V_{-\lm,\eps}$). We consider circle models of these representations where the group $K'_0$ acts by the standard rotations of $S^1$. We construct a (non-zero) explicit functional $l^{mod}_{E\otimes V_{-\tau,\eps'}}\in\Hom_G(E\otimes V_{-\tau,\eps'},\bc)$ which we call a model functional. Let us denote by
$sgn(\theta,\theta',\theta'')=sgn\left((\theta-\theta')(\theta-\theta'')(\theta'-\theta'')\right)$ the function taking values $\pm 1$ on $(S^1)^3\setminus \{$points with at least $2$ coordinates equal$\}$, and changes the sign at points removed. There are exactly two open orbits for the diagonal action of $\operatorname{PGL}^+_2(\br)$ on $(S^1)^3$, and the above function is an invariant of these. Moreover it is antisymmetric with respect to the action of the element $\dl$. It is easy to see (as in \cite{BR2})  that in the circle model of class one
representations the kernel of $l^{mod}_{E\otimes V_{-\tau,\eps'}}$ is
given by the following function in three variables $\theta,\
\theta',\ \theta''\in S^1$
\begin{eqnarray}\label{K-circle}
K_{\lm,-\lm,\tau}(\theta,\theta',\theta'')= (sgn(\theta,\theta',\theta''))^{\eps'}
\cdot&&\\|\sin(\theta-\theta')|^{\frac{-1-\tau}{2}}|\sin
(\theta-\theta'')|^{\frac{-1-2\lm+\tau}{2}}|\sin
(\theta'-\theta'')|^{\frac{-1+2\lm+\tau}{2}}\ .\nonumber
\end{eqnarray} (The factor should be equal to $(sgn(\theta,\theta',\theta''))^{\eps+\eps+\eps'}$, but it gives the same function.) This function also defines the kernel of the map $T_\tau=T^{mod}_{\tau,\eps'}:E\to V_{\tau,\eps'}$ via
the relation
$$\langle T_\tau(w),v\rangle_{V_{\tau,\eps'}}=\frac{1}{(2\pi)^3}\int_{(S^1)^3}w(\theta,\theta')v(\theta')
K_{\lm,-\lm,\tau}(\theta,\theta',\theta'')d\theta
d\theta'd\theta''\ . $$

\subsubsection{Integral transform} The model functional $d_\tau$ is given as the scalar product with a norm one $K'_0$-invariant vector $e'_0(\theta'')\equiv 1$ (i.e., $d_\tau(v)=\langle v, e'_0\rangle_{V_{\tau,\eps'}}$ as in section \ref{T-invar-funct-sect}).
Hence we have
\begin{eqnarray}\label{e'0}d_{\tau}(T_\tau(w))=\langle T_\tau(w),e'_0\rangle_{V_\tau}
=\\ \frac{1}{(2\pi)^3}\int
w(\theta,\theta')K_{\lm,-\lm,\tau}(\theta,\theta',\theta'')e'_0(\theta'')d\theta
d\theta'd\theta'', \nonumber\end{eqnarray}
for any $w\in C^\8(S^1\times S^1)$. For what follows, it will be enough to assume
that the vector $w\in E$ is $\Dl K'_0$-invariant. Such a
vector $w$ can be described by a function of one variable; namely,
$w(\theta,\theta')= u(c)$ for $u \in C^\8(S^1)$ and $c = (\theta -
\theta')/2$. We have then $\hat w(n,-n)=\hat
u(n)=\frac{1}{2\pi}\int_{S^1} u(c)e^{-inc}dc$ -- the Fourier
transform of $u$.

We introduce a new kernel (note that  $e'_0(\theta'')\equiv 1$ in \eqref{e'0})
\begin{eqnarray}\label{k-lm-def}
k_\tau(c)=k_{\lm,\tau} \left({\scriptstyle \frac{\theta -
\theta'}{2}}\right)=\frac{1}{2\pi}\int_{S^1}
K_{\lm,-\lm,\tau}(\theta,\theta',\theta'')d\theta''
\end{eqnarray} and the corresponding integral transform
\begin{eqnarray}u^\sharp(\tau)=u^\sharp_\lm(\tau)=\frac{1}{(2\pi)^2}\int_{S^1}
u(c)k_{\tau}(c)dc\ ,
\end{eqnarray} suppressing the dependence on $\lm$ as we have fixed
the representation $V_{\lm,\eps}$. The transform is clearly defined for
any smooth function $u\in C^\8(S^1)$, at least for $\tau\in i\br$.
In fact, it could be defined for all $\tau\in \bc$, by means of
analytic continuation. We will discuss this in \ref{DiscrSer} where we deal with discrete series.

Note that $k_\tau$ is the average of the kernel
$K_{\lm,-\lm,\tau}$ with respect to the action of $\Dl K'_0$, or, in
other terms, $k_\tau=T_\tau^*(e'_0)\in E^*$ is the pullback of the $K_0'$-invariant vector $e'_0\in
V_{\tau,\eps'}$ under the map $T_\tau^*$.

\subsection{Gelfand formation: spectral identity} We invoke now the uniqueness principle $\dim\Hom_G(E\otimes V_{-\tau,\eps'},\bc)\leq 1$. This implies that the two invariant functionals, $l^{aut}_{i}$ and $l^{mod}_{E\otimes V_{-\tau}}$, we constructed are proportional. Hence there are constants $\beta(i)=\beta(\nu_i)\in\bc$ such that
\begin{eqnarray}\label{aut-mod}
T^{aut}_{i}=\beta(i)\cdot T_{\tau_i}\ .
\end{eqnarray}
We also recall that the automorphic $K'_0$-invariant functional on $V_{\tau_i,\eps_i}$ satisfies the relation $d_i^{aut}=\alpha(i)\cdot d_{\tau_i}$ as in \eqref{T-invar-funct-sect}. Denoting by  $\g(i)=\alpha(i)\beta(i)$ the product of these constants, we
rewrite the spectral expansion \eqref{h-lm-spec} in the form
\begin{eqnarray}\label{h-lm-spec2}  h_\lm(w_u)=\sum_{\tau_i} d^{aut}_i(T^{aut}_i(w))=\sum_{\tau_i} \g(i)d_{\tau_i}(T_{\tau_i}(w_u))\\ =\sum_{\tau_i} \g(i)\cdot u^\sharp(\tau_i),\nonumber
\end{eqnarray} for any $\Dl K'_0$-invariant vector $w_u\in E$ corresponding to the function $u \in C^\8(S^1)$. This is the spectral expansion corresponding to the triple product.
Taking this together with the abelian spectral expansion \eqref{spec-h-lm}, we obtain the spectral identity corresponding to {\it two} Gelfand pairs (we call such a pair a Gelfand formation, and the corresponding identity the period identity):
\begin{eqnarray}\label{rs-K-u}
\sum_n|a_n|^2 \hat u(n)=  u(0)+
\sum_{\tau_i\not=1}\g(i)\cdot u^\sharp(\tau_i)\ .
\end{eqnarray}
Here on the right we singled out  the contribution from the
trivial representation (i.e., $\tau=1$), and wrote it in the form
$u(0)=(\Vol(\ck)/\Vol(X)^\haf )
\cdot  u(0)$ under our
normalization of measures $\Vol(X)=1$ and $\Vol(\ck)=1$.

By choosing the appropriate test function $u$, and  estimating the right hand side in \eqref{rs-K-u}, we will obtain the bound \eqref{short-bound-an}.

\subsection{Bounds for spherical Fourier coefficients} We are
interested in getting a bound for the coefficients $a_n$.
The idea of the proof of the crucial bound \eqref{short-bound-an} is to find a test vector $w\in V\otimes \bar
V$, i.e., a function $w\in C^\8(S^1\times S^1)$, such that when
substituted in the period identity \eqref{rs-K-u} it will
produce a weight $\hat w$ which is not too small for a given $n$,
$|n|\to\8$. We then have to estimate the spectral density of such
a vector, i.e., the transform $w^\sharp$. One might be tempted to
take $w$ such that $\hat w$ is essentially a delta function (i.e.,
picks up just a few coefficient $a_n$ in \eqref{rs-K-u}).
However, for such a vector we have no means to estimate the right
hand side of the  formula \eqref{rs-K-u} because $w^\sharp$ is
spread over a long interval of the spectrum. The solution to this
problem is well known in harmonic analysis. One takes a function
which produces a weighted sum of the coefficients $|a_k|^2$
for $k$ in a certain range depending on $n$ and such that its
transform $w^\sharp$ is spread over a shorter interval. For such
test vectors $w$, we give an essentially sharp bound for the value of the diagonal period
$d_{\Dl\mathcal{K}}(w)=\int_\ck r_\Dl(\nu_E(w))d\ck$.

\subsubsection{Proof of bound \eqref{short-bound-an}} \label{proof-thm-3}We
start with the period identity \eqref{rs-K-u} and construct
an appropriate $\Dl K_0'$-invariant vector $w\in E$, i.e., a function
$u\in C^\8(S^1)$ such that $w(\theta,
\theta')=u((\theta-\theta')/2)$.
We have the following elementary  technical lemma.

\begin{lemma}\label{lem2}For any
integers $N\geq T\geq 1$, there exists a smooth function
$u_{N,T}\in C^\8(S^1)$ such that
\begin{enumerate}
    \item $|u_{N,T}(0)|\leq \al T$,
    \item $\hat u_{N,T}(k)\geq 0$ for all $k$,
    \item $\hat u_{N,T}(k)\geq 1$ for all $k$ satisfying $|k-N|\leq T$,
    \item $|u_{N,T}^\sharp(\tau)|\leq \al
    T|N|^{-\haf}(1+|\tau|)^{-\haf}+\al T(1+|\tau|^{-5/2})$ for $|\tau|\leq N/T$,
    \item $|u_{N,T}^\sharp(\tau)|\leq \al
    T(1+|\tau|)^{-5/2}$ for $|\tau|\geq N/T$,
\end{enumerate} for some fixed constant $\al>0$ independent of $N$ and $T$.
\end{lemma}

The proof of this Lemma is given in Appendix A of \cite{2}. One constructs the
corresponding function $u_{N,T}$ by considering a function of the
type $u_{N,T}(c)=Te^{-iNc}\cdot\left(\psi*\bar\psi\right)(Tc)$
with a fixed smooth function $\psi\in C^\8(S^1)$ of a support in a
small interval containing $1\in S^1$ (here $*$ denotes the
convolution in $C^\8(S^1)$). Such a function obviously satisfies
conditions $(1)-(3)$ and the verification of $(4)-(5)$ is reduced
to a routine application of the stationary phase method (similar
to our computations in \cite{BR3}). \qed

We return to the proof of bound \eqref{short-bound-an}. In the proof we will use
two bounds for the coefficients $\alpha(i)$ and $\beta(i)$.
Namely, it was shown in \cite{BR2} that
\begin{eqnarray}\label{convex-in-proof}
\sum\limits_{A\leq|\tau_i|\leq 2A}|\beta(i)|^2\leq a A^2\ ,
\end{eqnarray}  for
any $A\geq 1$ and some explicit $a>0$.

The second  bound  that we will need is the bound
\begin{eqnarray}
\label{Horm}\sum_{A\leq|\tau_i|\leq 2A}|\alpha(i)|^2\leq b A^2\
,\end{eqnarray} valid for any $A\geq 1$ and some $b$. In disguise
this is the classical bound of L. H\"{o}rmander \cite{Ho} for the
average value at a point for eigenfunctions of the
Laplace-Beltrami operator on a compact Riemannian manifold (e.g.,
$\Dl$ on $Y$). This follows from the normalization
$|\beta(i)|^2=|\phi'_{\tau_i}(x_0)|^2$ we have chosen in
\eqref{b-0-nomalize} for $K'_0$-invariant eigenfunctions. In fact,
the bound \eqref{Horm} is standard in the theory of the Selberg
trace formula (see \cite{Iw}), and also can be easily deduced from
considerations of \cite{BR3}.

We plug a test function satisfying conditions $(1)-(5)$ of Lemma
\ref{lem2} into the identity \eqref{rs-K-u}. Using
the Cauchy-Schwartz inequality and taking into account bounds
\eqref{convex-in-proof} and \eqref{Horm}, we obtain

$\begin{array}{lr}
 \sum\limits_{|k-N|\leq T}|a_k|^2\leq \sum\limits_k|a_k|^2 \hat u_{N,T}(k)
= u_{N,T}(0) +\sum\limits_{\tau_i\not=1}\alpha(i) \beta(i) u_{N,T}^\sharp(\tau_i) &\\
\leq  \al T+\sum\limits_{|\tau_i|\leq N/T}\al T|N|^{-\haf}(1+|\tau_i|)^{-\haf}|\alpha(i) \beta(i)|
+\sum\limits_{\tau_i\not=1}\al T(1+|\tau_i|)^{-5/2}|\alpha(i) \beta(i)|&\\
\leq \al T+\al T|N|^{-\haf}\sum\limits_{|\tau_i|\leq N/T}(1+|\tau_i|)^{-\haf}\left(|\alpha(i)|^2+|\beta(i)|^2\right)
   +&\\\al T\sum\limits_{\tau_i\not=1}(1+|\tau_i|)^{-5/2}\left(|\alpha(i)|^2+
  |\beta(i)|^2\right)
\leq  \al T+CT|N|^{-\haf}(N/T)^{3/2+\eps}+DT=\\
  c'T+CT^{-\haf-\eps}|N|^{1+\eps}\ ,
\end{array}$
\\ for any $\eps>0$ and some constants $c',\ C,\ D>0$.

Setting $T=N^{2/3}$, we obtain $ \sum_{|k-N|\leq
N^{2/3}}|a_k|^2\leq A_\eps N^{2/3+\eps}$ for any $\eps>0$.

This finishes the proof of the bound \eqref{short-bound-an}, and with it the proof of Theorem \ref{claimB} for representations of principal series. In \ref{DiscrSer} we explain the case of discrete series.
\qed

\appendix
\section{Representations of the discrete series}
\label{DiscrSer}
The aim of this appendix is to indicate changes needed in order to treat discrete series representations. The main difference in treatment of the discrete  series is the lack of their convenient realization. The problem is computational and not conceptual. We show how to reduce necessary computations to the case of induced representations. With these changes, the treatment is identical to the case of induced representations we discussed in the main body of the paper. The main reason why we can make this reduction comparatively easy is that we assume that the representation is fixed. Similar computations that take into account the weight of the representation would be  more complicated.

\subsection{Representations and their realizations} Let $k\geq 2$ be an even integer, and $(D_k,\pi_k)$ be the corresponding discrete series representation of $\PGLR$. In particular, for $m\in 2\bz$, the space of  $K_0$-types of weight $m$ is non-zero (and in this case is one-dimensional) if and  only if $|m|\geq k$. This defines $\pi_k$ uniquely. Under the restriction to $\PSLR$, the representation $\pi_k$ splits into two representations $(D^\pm_k,\pi^\pm_k)$ of ``holomorphic" and ``anti-holomorphic" discrete series  of $\PSLR$, and  the element $\dl$ interchanges these.

We consider two realizations of discrete series as subrepresentations and as quotients of induced representations. Consider the space  $\CH_{k-2}$  of smooth even homogeneous functions on $\br^2\setminus 0$ of the homogeneous degree $k-2$ (i.e., $f(tv)=t^{k-2}f(v)$ for any $t\in \br^\times$ and $0\not=v\in\br^2$). We have the natural action of $\GLR$ given by $\tilde{\pi}_{k}(g)f(v)=f(g\inv v)\cdot \det(g)^{(k-2)/2}$, which is trivial on the center and hence defines a representation of $\PGLR$. There exists the unique non-trivial invariant subspace $V_{k-2}\subset \CH_{k-2}$. The space $V_{k-2}$ is finite-dimensional, $\dim V_{k-2}=k-1$, and is generated by the homogeneous polynomials of degree $k-2$. The quotient space $\CH_{k-2}/ V_{k-2}$ is isomorphic to the space of smooth vectors of the discrete series representation $\pi_k$. Note that $V_{k-2}$ consists of vectors with $K_0$-type in the range $|n|<k$.

We also consider the dual situation. Let  $\CH_{-k}$ be the space of smooth even homogeneous functions on $\br^2\setminus 0$ of the homogeneous degree $-k$. There is a natural $\PGLR$-invariant pairing $\langle\ ,\ \rangle: \CH_{k-2}\otimes \CH_{-k}\to\bc$. Hence $\CH_{-k}$ is the smooth dual of $\CH_{k-2}$, and vice versa. There exists the unique  non-trivial invariant subspace in $\CH_{-k}$, and it is isomorphic to $D_k$. The quotient $\CH_{-k}/D_k$ is isomorphic to the finite-dimensional representation $V_{k-2}$.

Consider the restriction of smooth functions in $\CH_{k-2}$ to functions on the circle $S^1\subset \br^2\setminus 0$. There exists a unique (up to a multiple) Hermitian $G$-invariant form on $\CH_{k-2}\simeq C^\8_{ev}(S^1)$ given by $F_k(f,g)=\sum\limits_{n=-\8}^{n=-k}\g_-(k,n)\cdot a_n\bar b_n+\sum\limits^{n=\8}_{n=k}\g_+(k,n)\cdot a_n\bar b_n$, where $a_n$ and $b_n$ are the Fourier coefficients of $f$ and $g$ correspondingly, and the coefficients $\g_\pm(k,n)$ are given  by $\g_\pm(k,n)=\G(k\pm n)/\G(\pm n-k)$. (Note that $F_k$ is degenerate on  $\CH_{k-2}$.) The Hermitian form $F_k$ corresponds to the invariant Hermitian product on $\pi_k$ after passing to the quotient (in particular, it is positive definite on the quotient). We denote it $\langle \ ,\ \rangle_{D_k}$. Note that asymptotically the coefficients $\g_\pm(k,n)$ grow as $|n|^{2k+1}$, and hence the corresponding norm on $\pi_k$ resembles  the Sobolev norm on $S^1$. We will not use the invariant norm, and instead use the natural {\it linear} pairing on  $D_k\otimes D^*_k$ and its relation to the natural pairing on $\CH_{k-2}\otimes \CH_{-k}$.

Consider the representation  $D_k$ and its (smooth) dual $D_k^*$ (of course $D_k\simeq D^*_k$, but we want to distinguish between these two copies). We want to find a way to compute the canonical  pairing $\langle\ ,\ \rangle_k:
D_k\otimes D^*_k\to\bc$ in terms of the canonical pairing $\langle\ ,\ \rangle:\CH_{k-2}\otimes \CH_{-k}\to\bc$. Consider  the imbedding $i_k: D_k^*\to\CH_{-k}$ and the quotient map $q_k: \CH_{k-2}\to D_k$. Hence we also have
$id\otimes i_k: \CH_{k-2}\otimes D_k^*\to\CH_{k-2}\otimes\CH_{-k}$ and $q_k\otimes id: \CH_{k-2}\otimes\CH_{-k}\to D_k\otimes \CH_{-k}$. The image of the composition $ q_k\otimes id\circ id\otimes i_k$ is equal to the image of the imbedding $id\otimes i_k: D_k\otimes D^*_k\to D_k\otimes \CH_{-k}$.  Hence for a pair of vectors $u\otimes v\in \CH_{k-2}\otimes\CH_{-k}$ such that $v=i_k(v')$, $v'\in D^*_k$,  we have $\langle q_k(u), v\rangle_k=\langle u,v'\rangle$.  This will be our way to compute the invariant Hermitian norm on $D_k$. The important difference, as opposed to the invariant Hermitian form $\langle \ ,\ \rangle_{D_k}$ on $D_k$,   is that it is not possible to say what is the norm of a vector $v\in D_k$ (for this, one uses the intertwining operator between $\CH_{-k}$ and $\CH_{k-2}$ which we want to avoid). It is however possible to say when the scalar product between a vector and a co-vector is one. In general this does not allow one to construct a projector to a particular vector. However in our situation we will need projectors onto pure $K_0$-types, and since these are essentially unique the above construction is enough for our purposes. We note, however, that for a general pair $u\otimes v \in \CH_{k-2}\otimes\CH_{-k}$ not satisfying the condition $v=i_k(v')$ (i.e., $v\not\in Im(i_k)\subset \CH_{-k}$), the value of $\langle u,v\rangle$ might be far from a pairing of vectors in the discrete series (e.g., paring of high $K_0$-types components of $u,\ v$). This is because components belonging  to the finite dimensional representations might be dominant in the pairing $\langle\ ,\ \rangle$.

We now compare this to the automorphic picture. Let $\nu: D_k\to L^2(X)$ be an isometric automorphic realization of a discrete series representation. Consider the complex conjugate realization $\bar\nu:D^*_k\to L^2(X)$. Integration along the diagonal $\Delta X\hookrightarrow X\times X$ corresponds to the canonical automorphic pairing  $\langle\ ,\ \rangle_k^{aut}:
D_k\otimes D^*_k\to\bc$.  We have  $\langle\ ,\ \rangle_k^{aut}= \langle\ ,\ \rangle_k$  since it is assumed that $\nu$ is an isometry. While as we noted it does not allow us to compute norm of a vector in $D_k$, we can produce some functions on $X$ corresponding to norm one vectors. For example, let $e_n\subset D_k$ be a vector of  $K_0$-type $n$, and $e^*_{-n}\in D^*_k$ the dual $K_0$-type vector (i.e., $\langle e_n,e^*_{-n}\rangle_k=1$). We have then $\nu(e_n)(x)\bar\nu(e^*_{-n})(x)=|\nu(\tilde e_n)(x)|^2$ as functions on $X$,  where $\tilde e_n\in D_k$ is a norm one $K_0$-type.

\subsection{Invariant functionals} We now construct $K_0$-equivariant functionals and triple functionals for discrete series using the space $\CH_{k-2}\otimes\CH_{-k}$.
Let $\nu\otimes\bar\nu: D_k\otimes D_k^*\to L^2(X\times X)$ be as before, and $\ck\subset X$ be a $K'_0$-orbit. For a character $\chi:K'_0\to S^1$, we define the automorphic functional $d_{\chi_n}^{aut}(\nu(u))=\int_\ck\nu(u)(k)\bar\chi_n(k)d\ck$, which is proportional to the model functional $d_{n}(u)=\langle u, e'_n\rangle_{D_k}$, i.e., $d_{\chi_n}^{aut}=a_n\cdot d_n$, where  $e'_n\in D_k$ is a norm one vector of $K'_0$-type $n$. On the other hand, we have the corresponding functional  $\hat d_{n}\otimes \hat d_{-n}:\CH_{k-2}\otimes\CH_{-k}\to \bc$ given by $\hat d_{n}\otimes \hat d_{-n}(u\otimes w)=\langle u, \hat e'_n\rangle\cdot \langle w, \hat e'_{-n}\rangle $, where   $\hat e'_n\otimes \hat e'_{-n}\in \CH_{-k}\otimes\CH_{k-2}$  is any pair of corresponding $K'_0$-types which are in duality $\langle  \hat e'_n, \hat e'_{-n}\rangle=1$. The main point here is that we can easily compute the value of $\hat d_{n}\otimes \hat d_{-n}(u\otimes w)$ since in the circle model $\CH_{-k}\otimes\CH_{k-2}\simeq C^\8(S^1\times S^1)$ these are Fourier coefficients of functions $u$ and $w$.   For $|n|\geq k$, we have $d_{\chi_n}^{aut}\otimes d_{\chi_{-n}}^{aut}=|a_n|^2\cdot\hat d_{n}\otimes \hat d_{-n}$.
Let $e_j\in D_k$ be a norm one $j$-th $K_0$-type (we recall that there are two {\it different} maximal compact connected subgroups $K_0=\PSOt$ and $K_0'$ involved in the description of the set $O\subset X$). As in \eqref{specH-lm-aut}, we are interested in coefficients $|d_n(e_j)|^2=|\langle e_j, e'_n\rangle_{D_k}|^2$. For $|j|\geq k$, a vector with the $K_0$-type equal to $j$ automatically belongs to the subspace $D_k\subset\CH_{-k}$, and we have $|d_n(e_j)|^2=d_n(e_j)\cdot d_{-n}(e_{-j})=\hat d_{n}(\tilde e_j)\cdot \hat d_{-n}(\tilde e_{-j})$.  Here $e_{\pm j}\in D_k$ are norm one vectors and $\tilde e_{\pm j}\in \CH$ are any vectors in duality $\langle  \tilde e_j, \tilde e'_{-j}\rangle=1$. This implies that we have the integral representation for  the quantity similar to the integral \eqref{osc-int} we used for principal series. In particular, this function exhibits the same Airy type behavior. Hence as for the principal series, we need  to estimate the coefficients $a_n$. For this we construct model triple products as in the case of the principal series.

Let $(V_{\tau,\eps'},\pi_{\tau,\eps'})$ be a unitary representation of the principal series. The space $\Hom_G(D_k\otimes D^*_k, V_{\tau,\eps'})$ is one-dimensional. We will work with the space of invariant trilinear functionals $\Hom_G(D_k\otimes D^*_k\otimes V_{-\tau,\eps'},\bc)$ instead. We first construct explicitly a non-zero element $\hat l_{k,\tau,\eps}$ in the space $\Hom_G(\CH_{k-2}\otimes\CH_{-k}\otimes V_{-\tau,\eps'},\bc)$ (which is in fact also is one-dimensional), and then use it to define a non-zero element in the space  $\Hom_G(D_k\otimes D^*_k\otimes V_{-\tau,\eps'},\bc)$. What is more important, we will use $\hat l_{k,\tau,\eps}$ to carry out our computations in a similar way to the principal series.

Consider the following function (compare to \eqref{K-circle}) in three variables $\theta,\
\theta',\ \theta''\in S^1$
\begin{eqnarray*}\label{K-circle-discr}
K_{k-2,-k,\tau}(\theta,\theta',\theta'')=(sgn(\theta,\theta',\theta''))^{\eps'}
\cdot\\ |\sin(\theta-\theta')|^{\frac{-1-\tau}{2}}|\sin
(\theta-\theta'')|^{\frac{-1+\tau}{2}-k+1}|\sin
(\theta'-\theta'')|^{\frac{-1+\tau}{2}+k-1}\ .
\end{eqnarray*}  Viewed as a kernel, it defines an invariant non-zero functional on the (smooth part of) the representation $\CH_{k-2}\otimes\CH_{-k}\otimes V_{-\tau,\eps'}\simeq C^\8_{ev}(S^1\times S^1\times S^1)$.
 Such a kernel should be understood in the regularized sense (e.g.,  analytically  continued following  \cite{G1}). However, since we are interested in $\tau\in i\br$, all the exponents in \eqref{K-circle-discr} are non integer. This implies that the regularized kernel coincides with the above function when integrated against any test function vanishing in a neighborhood of singularities of the kernel \eqref{K-circle-discr}.  We denote the corresponding functional by $\hat l_{k,\tau,\eps}\in \Hom_G(\CH_{k-2}\otimes\CH_{-k}\otimes V_{-\tau,\eps'},\bc)$. Such a functional defines a non-zero functional $l_{k,\tau,\eps}\in \Hom_G(D_k\otimes D^*_k\otimes V_{-\tau,\eps'},\bc)$ by requiring that $l_{k,\tau,\eps}$ and $\hat l_{k,\tau,\eps}$  induce the same functional on the space $\CH_{k-2}\otimes D^*_k\otimes V_{-\tau,\eps'}$.   We call $l_{k,\tau,\eps}$ the model functional for the discrete series. The difference with principal series clearly lies in the fact that we only can compute the auxiliary functional $\hat l_{k,\tau,\eps}$. However, for $k$ fixed, it turns out that necessary computations are identical to the ones we performed for the principal series.

 As  in \eqref{k-lm-def}, we are interested in
$k_{k,\tau} \left({\scriptstyle \frac{\theta -
\theta'}{2}}\right)=\frac{1}{2\pi}\int_{S^1}
K_{k-2,-k,\tau}(\theta,\theta',\theta'')d\theta''\ .$
Moreover, we are interested in asymptotic of this function as $|\tau|\to\8$ and $k$ is {\it fixed}. The analysis of such an integral is standard in the theory of the stationary phase method, and is similar (since $k$ is fixed) to one which is done  in \cite{2} and \cite{BR3}. This is because the contribution from any small enough neighborhood of the singularity of the above kernel is negligible when $|\tau|\to\8$ due to the high oscillation of the kernel. The situation is of course different if also $k\to\8$, but this is not relevant to our present situation.   Hence we obtain the asymptotic expansion of the trilinear invariant functional essentially identical to the one we obtained for principal series. This allows us to use the same calculations which we already explained before. In particular, we have to construct  test vectors analogous to the ones obtained from  functions $u_{N,T}\in C^\8(S^1)$ in  Lemma \ref{lem2}. We can take the same function as for the principal series, but while it gives a $\Dl K_0$-invariant vector in the space $\CH_{k-2}\otimes\CH_{-k}$ such a vector does not belong to $\CH_{k-2}\otimes D^*_k$. As a result, it is difficult to estimate the norm of such a vector. However it is easy to correct the function $u_{N,T}$ so that it will have vanishing Fourier coefficients $\hat u_{N,T}(n)$ for $|n|< k$ (and hence will produce a vector in $\CH_{k-2}\otimes D^*_k$). Since we are interested in the asymptotic as $N,\ T\to\8$,  this finite correction does not affect it. Here  the fact that we consider {\it fixed} $k$ is crucial. Without this restriction, it is more difficult to calculate the corresponding value of the trilinear invariant functional (an in fact the resulting asymptotic should be different).

\begin{acknowledgements} It is a pleasure to thank Joseph Bernstein for endless discussions concerning automorphic functions, Leonid  Polterovich for turning my attention to the paper \cite{Ta} and illuminating discussions, David Kazhdan who suggested to put the results  of the paper in writing, and  Elon Lindenstrauss for useful remarks.
\end{acknowledgements}

\end{document}